\documentclass[a4paper,10pt]{article}

\usepackage[T1]{fontenc}
\usepackage[latin1]{inputenc}
\usepackage[english]{babel}
\usepackage{amsmath,amsthm,amssymb,geometry,fancyhdr,stackrel}
\usepackage[all]{xy}

\title{Topological invariants of piecewise hereditary algebras}
\author{Patrick Le Meur}
\date{\today}

\theoremstyle{definition}
\newtheorem{definition}{Definition}[section]
\theoremstyle{plain}
\newtheorem{Thm}[]{Theorem}
\newtheorem{Cor}[Thm]{Corollary}
\newtheorem{cor}[definition]{Corollary}
\newtheorem{lem}[definition]{Lemma}
\newtheorem{prop}[definition]{Proposition}
\theoremstyle{remark}

\newtheorem{rem}[definition]{Remark}

\def\ts#1{\normalfont{\textsf{#1}}}
\def\w#1{\widetilde{#1}}
\def\c#1{\mathcal{#1}}
\def\dbh{\c D^b(\c H)}
\def\o#1{\overline{#1}}
\def\l{F_{\lambda}}
\def\sq{\null\hfill$\square$\\}
\def\dbc{\c D^b(\ts{mod}\, \c C)}
\def\dba{\c D^b(\ts{mod}\, A)}
\def\dbb{\c D^b(\ts{mod}\, B)}

\makeatletter
\renewcommand\section{\@startsection {section}{1}{0mm}%
{-3.5ex \@plus -1ex \@minus -.2ex}
                                   {0.5ex \@plus.2ex}%
                                   {\normalfont\large\bfseries}}
\makeatother
\begin{document}
\sffamily
\maketitle
\abstract{
We investigate the Galois coverings of
piecewise algebras and more
particularly their behaviour under derived
equivalences. Under a technical assumption which is satisfied if the
algebra is derived equivalent to a hereditary algebra, we prove that there
exists a universal Galois covering whose group of automorphisms is
free and depends only on the derived category of the algebra. As a
corollary, we prove that the algebra is simply connected if and only
if its first Hochschild cohomology vanishes.}

\section*{Introduction}$\ $

Let $k$ be an algebraically closed field and $A$ a basic
finite dimensional $k$-algebra (or, simply, an algebra). The representation theory
studies the
category $\ts{mod}\, A$ of finite dimensional (right) $A$-modules and
also its bounded derived category $\dba$. From this point of view,
some classes of algebras play an important rôle: The
\textit{hereditary} algebras, that is, path algebras $kQ$ of finite quivers $Q$
with no oriented cycle; the \textit{tilted} algebras, that is, of the
form $\ts{End}_{kQ}(T)$, where $T$ is a tilting $kQ$-module; and, more generally, the
\textit{piecewise hereditary} algebras, that is, the algebras $A$ such
that $\dba$ is triangle equivalent to $\dbh$ where
$\c H$ is a $\ts{Hom}$-finite hereditary abelian category with split
idempotents (if
$\c H=\ts{mod}\, kQ$ then $A$ is called  piecewise hereditary \textit{of type $Q$}). These
algebras are
particularly well understood (see
\cite{assem_simson_skowronski,happel_book,simson_skowronski,simson_skowronski2},
for instance).

The piecewise hereditary algebras arise in many areas of representation theory. For
example,
the cluster category $C_A$ of a piecewise hereditary algebra was
introduced in \cite{bmrrt,ccs} as a tool to study conjectures related to cluster algebras (\cite{fz}). 
Another example is the study of self-injective
algebras, that is, algebras $A$ such that $A\simeq DA$ as right $A$-modules. Indeed,
to any algebra $A$ is associated the
repetitive category $\widehat{A}$, which is a Galois covering with
group $\mathbb{Z}$ of the trivial extension $A\ltimes DA$ (see \cite{hw}). Assume that some group $G$ acts
freely on $\widehat{A}$ thus defining a Galois covering $\widehat{A}\to \widehat{A}/G$
with group $G$. If $\widehat{A}/G$ is a finite dimensional
algebra, that is, if it has finitely many objects as a category, then
it is self-injective and called \emph{of type $Q$} if $A$ is tilted of
type $Q$. It is proved in
\cite{skowronski} that any self-injective algebra of polynomial growth
and admitting a Galois covering by a strongly simply connected category is of
the form $\widehat{A}/G$ for some tilted algebra $A$ and some
infinite cyclic group $G$. The class of self-injective algebras of type
$Q$ has been the object of many studies recently (see \cite{sy,sy2,sy3}).

In this text we investigate the Galois coverings of piecewise hereditary
algebras. The Galois coverings of algebras and, more generally, of
$k$-categories, were introduced in
\cite{gabriel,riedtmann} for the classification of
representation-finite algebras. Consider $A$ as a locally bounded $k$-category: If
$1=e_1+\ldots+e_n$ is a decomposition of the unity into primitive
orthogonal idempotents, then $\ts{ob}(A)=\{e_1,\ldots,e_n\}$ and the space
of morphisms from $e_i$ to $e_j$ is $e_jAe_i$.
Then a \emph{Galois covering} of the $k$-category $A$ is a $k$-linear
functor $F\colon\c C\to A$ where $\c C$ is a $k$-category endowed with
a free action of $G$, that is, $G$ acts freely on $\ts{ob}(\c C)$, such
that $F\circ g=F$ for every $g\in G$ and the induced functor $\c
C/G\to A$ is an isomorphism (\cite{gabriel}). In such a situation,
$\ts{mod}\,\c C$ and $\ts{mod}\, A$ are related by the so-called push-down functor
$F_{\lambda}\colon \ts{mod}\, \c C\to \ts{mod}\, A$, that is, the
extension-of-scalars functor. Often,
$F_{\lambda}$ allows nice comparisons between $\ts{mod}\, \c C$ and
$\ts{mod}\, A$. For example: The action of $G$ on $\c C$ naturally defines
an action $(g,X)\mapsto \,^gX$ of $G$ on $\c C$-modules. When this action is free on
indecomposable $\c C$-modules, $F_{\lambda}$ defines an isomorphism of
translation quivers
between $\Gamma(\ts{mod}\, \c C)/G$ and a union of some components of the
Auslander-Reiten quiver $\Gamma(\ts{mod}\, A)$ of $A$
(see \cite{dowbor_skowronski,gabriel}). 

The comparisons allowed by the
covering techniques raise naturally the following questions: Given an
algebra $A$, is it possible to describe all the Galois coverings of
$A$ (in particular, does $A$ admit a universal Galois covering, as
happens in topology)? Is it
possible to characterise the simple connectedness of $A$
(that is, the fact that $A$ has no proper Galois covering by a
connected and locally bounded $k$-category)? In view of the above
discussion on self-injective algebras, these questions are
particularly relevant when $A$ is piecewise hereditary of type $Q$. In case
$A=kQ$, the answers are well-known: The Galois
coverings of $kQ$ correspond to the ones of the underlying graph of
$Q$; and $kQ$ is simply connected if and only if $Q$ is a tree, which
is also equivalent to the vanishing of the first Hochschild cohomology
group $\ts{HH}^1(kQ)$ (\cite{cibils_h1}). Keeping in mind the general
objective of representation theory,
one can wonder if the data of the Galois coverings of $A$ and the
simple connectedness of $A$ depend only on the
bounded derived category $\dba$. Again, it is natural to treat
this problem for piecewise hereditary algebras.
Up to now, there are no general solutions to the above problems. The
question of the description of the Galois coverings and the one of
the characterisation of simple connectedness have found a satisfactory
answer in the case of standard representation-finite algebras
(see \cite{buchweitz_liu,gabriel}). This is mainly due to the fact that
the Auslander-Reiten quiver is connected and completely describes the
module category in this case. However, the infinite-representation case seems to be more
complicated. As an example, there exist string algebras
which admit no universal Galois covering (\cite{lemeur3}).
In the present text, we study the above problems when $A$ is piecewise
hereditary. As a main result, we prove the following theorem.
\begin{Thm}
\label{thm0.2}
Let $A$ be a connected algebra derived equivalent to a hereditary
abelian category $\c H$ whose oriented graph $\overrightarrow{\c K}_{\c H}$ of tilting
objects is connected. Then 
$A$ admits a universal Galois covering $\widehat{\c C}\to A$ with
group a free group $\pi_1(A)$ uniquely
determined by $\dba$. This
means that $\widehat{\c C}$ is connected and locally bounded and for
any Galois covering $\c C\to A$ with group $G$ where $\c C$ is
connected and locally bounded there exists a commutative diagram:
\begin{equation}
  \xymatrix@R=10pt{
\widehat{\c C} \ar@{->}[rd] \ar@{->}[dd]& &\\
&\c C \ar@{->}[d]&\\
A \ar@{->}[r]^{\sim} & A 
}\notag
\end{equation}
where the bottom horizontal arrow is an isomorphism extending the
identity map on $\ts{ob}(A)$. 
Moreover, $\widehat{\c
  C}\to\c C$ is Galois with group $N$ such that there is an exact
sequence of groups $1\to N\to \pi_1(A)\to G\to 1$.

Finally, if $A$ is hereditary of type $Q$ then $\pi_1(A)$ is the
fundamental group $\pi_1(Q)$ of the underlying graph of $Q$ and,
otherwise, the rank of $\pi_1(A)$ equals $\ts{dim}_k\,\ts{HH}^1(A)$
(which is $0$ or $1$).
\end{Thm}

We refer the reader to the next section for the definition of
$\overrightarrow{\c 
K}_{\c H}$. Recall (\cite{happel_unger3}) that the asumption on $A$ is satisfied
if $A$ is piecewise hereditary of type $Q$.

The above theorem implies that the Galois coverings of a piecewise
hereditary algebra are determined by the factor groups of $\pi_1(Q)$.
Also it shows that the
data of the Galois coverings is an invariant of the derived
category. Therefore so does the simple connectedness. Using the fact that the Hochschild
cohomology is invariant under derived equivalences
(see \cite{keller}), we deduce the following corollary of our main result.
\begin{Cor}
\label{cor0.1}
Let $A$ be as in Theorem~\ref{thm0.2}. The following
are equivalent:
\begin{enumerate}
\item[(a)] $A$ is simply connected.
\item[(b)] $\ts{HH}^1(A)=0$.
\end{enumerate}
If $A$ is piecewise hereditary of type $Q$, then (a) and (b) are also
equivalent to:
\begin{enumerate}
\item[(c)] $Q$ is a tree.
\end{enumerate}
\end{Cor}

 This corollary generalises some of the
results of \cite{assem_marcos_delapena,assem_skowronski} which studied
the same characterisation for tilted algebras of euclidean type and for
tame tilted algebras. Also, it gives a new class of
algebras for which the following question of Skowro\'nski
(\cite[Pb. 1]{skowronski2}) has a positive answer: \emph{Is $A$ simply
  connected if and only if $\ts{HH}^1(A)=0$?} Originally, this question was
asked for tame triangular algebras.

The methods we use to prove Theorem~\ref{thm0.2} allow us to prove the
last main result of this text. It shows that the Galois coverings have
a nice behaviour for piecewise hereditary algebras.
\begin{Thm}
\label{thm0.3}
Let $A$ be piecewise hereditary of
type $Q$ and  $F\colon \c C\to A$ be a Galois
covering with group $G$ where $\c C$ is connected and locally
bounded. Then $\c C$ is piecewise hereditary of type
a quiver $Q'$ such that there exists a Galois covering of quivers $Q'\to Q$ with
group $G$.
\end{Thm}

We now give some explanations on the proof of
Theorem~\ref{thm0.2}. For unexplained notions, we refer the reader to
the next section. 
Assume that $A$ is piecewise hereditary.
It is known from \cite[Thm. 2.6]{happel_reiten} that there exists an
algebra $B$ such that $A\simeq \ts{End}_{\dbb}(X)$ for some tilting
complex $X\in\dbb$ and such that $B$ has one of the following forms:
\begin{enumerate}
\item $B=kQ$, with $Q$ a finite quiver with no oriented cycle.
\item $B$ is a squid algebra.
\end{enumerate}
It is easy to check that Theorem~\ref{thm0.2} holds true for path algebras
of quivers and for squid algebras. Therefore 
 we are reduced to proving
that  Theorem~\ref{thm0.2} holds true
for $A$ and only if it holds true for $\ts{End}_{\dba}(T)$ for any tilting complex
$T\in\dba$. Roughly speaking, we
need a correspondence between the Galois coverings of $A$
and those of $\ts{End}_{\dba}(T)$.
 Therefore we use
 a construction introduced in \cite{lemeur5} for tilting modules:
 Given a Galois covering $F\colon\c 
C\to A$ with group $G$, the push-down functor $\l\colon\ts{mod}\,\c C\to
\ts{mod}\, A$ is exact and therefore induces an exact functor
$\l\colon\dbc\to\dba$. Also, the $G$-action on modules extends to a
$G$-action on $\dbc$ by triangle automorphisms. Now, let
$T\in\dba$ be a tilting complex and
$T=T_1\oplus\ldots\oplus T_n$ be an indecomposable
decomposition. Assume that the following conditions hold true for
every $i\in\{1,\ldots,n\}$:
\begin{enumerate}
\item[($H_1$)] There exists an indecomposable $\c C$-module
  $\widetilde{T}_i$ such that $F_{\lambda}\widetilde{T}_i=T_i$.
\item[($H_2$)] The stabiliser $\{g\in G\ |\
  ^g\widetilde{T}_i\simeq \widetilde{T}_i\}$  is the trivial group.
\end{enumerate}
Under these assumptions, the complexes $^g\w T_i$ (for $g\in G$ and
$i\in\{1,\ldots,n\}$) form a full subcategory of $\dbc$ which we
denote by $\ts{End}_{\dbc}(\w T)$. Then $\l\colon\dbc\to\dba$ induces a
Galois covering with group $G$:
\begin{equation}
\begin{array}{rclc}
\ts{End}_{\dbc}(\w T)&\to&
\ts{End}_{\dba}(T)&\\
^g\w T_i & \mapsto &T_i&\\
^g\w T_i\xrightarrow{u}\,^h\w T_j & \mapsto & T_i\xrightarrow{\l
  (u)}T_j& .
\end{array}\notag
\end{equation}
 Hence 
($H_1$) and ($H_2$) are technical conditions which allow one to
associate a Galois covering of $\ts{End}_{\dba}(T)$ to a Galois covering of
$A$. In particular, if $A$ admits a universal Galois covering, then
the associated Galois covering of $\ts{End}_{\dba}(T)$ is a good
candidate for being a universal Galois covering. This is indeed the
case provided that the following technical condition is satisfied:
\begin{enumerate}
\item[($H_3$)] If $\psi\colon
A\xrightarrow{\sim} A$ is an automorphism such that $\psi(x)=x$ for
every $x\in \ts{ob}(A)$, then $\psi_{\lambda}T_i\simeq T_i$, for every $i$.
\end{enumerate}
We therefore need to prove the assertions ($H_1$), ($H_2$)
and ($H_3$) for every Galois covering $F\colon\c C\to A$ and
every tilting complex $T\in\dba$.

The text is therefore organised as follows. In Section~\ref{sec:s1},
we recall some useful definitions and fix some notations. In
Section~\ref{sec:s2}, we define the exact functor
$\l\colon\dbc\to\dba$ associated to a Galois covering $F\colon\c C\to
A$. In
Section~\ref{sec:s3}, we introduce elementary transformations on
tilting complexes using approximations. The main result of the section
asserts that for $A$ piecewise hereditary of type $Q$ and  for any
tilting complexes $T,T'$, there exists a sequence
of elementary transformations relating $T$ and $T'$. We prove the
assertions ($H_1$), ($H_2$) and ($H_3$) in Section~\ref{sec:s4} using
the elementary transformations. We prove Theorem~\ref{thm0.3} as an
application of these results. Then, in Section~\ref{sec:s5}, we
establish a correspondence between the
Galois coverings of $A$ and those of $\ts{End}_{\dba}(T)$  for every
tilting complex $T$. Finally, we prove
Theorem~\ref{thm0.2} and Corollary~\ref{cor0.1} in Section~\ref{sec:s6}.

\section{Definitions and notations}
\label{sec:s1}
\paragraph{Modules over $k$-categories}$\ $

We refer the reader to \cite{bongartz_gabriel} for the definition of
$k$-categories and locally bounded $k$-categories. All locally bounded
$k$-categories
are assumed to be small and all functors between
$k$-categories are assumed to be $k$-linear (our module categories and
derived categories will be skeletally small).
Let  $\c C$ be a
$k$-category.
Following \cite{bongartz_gabriel}, a (right) $\c C$-module is a $k$-linear
functor  $M\colon\c C^{op}\to \ts{MOD}\, k$
where $\ts{MOD}\, k$ is the category of
$k$-vector spaces. The category of $\c C$-modules is denoted by $\ts{MOD}\,
\c C$. 
A module $M\in\ts{MOD}\, \c C$ is called \textit{finite
 dimensional} if $\sum\limits_{x\in \ts{ob}(\c
   C)}\ts{dim}_k M(x)<\infty$. The category of finite dimensional  $\c
 C$-modules is denoted by $\ts{mod}\, \c C$.
Note that the indecomposable projective $\c C$-module associated to
$x\in\ts{ob}(\c C)$ is the representable functor $\c C(-,x)$.
The projective dimension of a $\c C$-module $X$ is
denoted by $\ts{pd}_{\c C}(X)$. If
 $X\in\ts{mod}\, \c C$, then $\ts{add}(X)$ denotes the smallest full subcategory of
 $\ts{mod}\, \c C$ closed under direct summands and direct sums.
We refer the reader to \cite{assem_simson_skowronski} for notions on
tilting theory. If $A$ is an algebra, an $A$-module $T$ is
called \textit{tilting} if: (a) $T$ is multiplicity-free; (b)
$\ts{pd}_A(T)\leqslant 1$; (c)
$\ts{Ext}_A^1(T,T)=0$; (d) for every indecomposable projective $A$-module
$P$ there is an exact sequence $0\to P\to X\to Y\to 0$ in $\ts{mod}\, A$ where $X,Y\in
\ts{add}(T)$. Let $\c H$ be a hereditary abelian category. An object
$T\in\c H$ is called \emph{tilting} (see \cite{happel_reiten}) if: (a) $T$ is multiplicity-free;
(b) $\ts{Ext}_{\c H}^1(T,T)=0$; (c) whenever $\ts{Hom}_{\c
  H}(T,X)=\ts{Ext}^1_{\c H}(T,X)=0$ for $X\in\c H$, then $X=0$. The
set of isomorphism classes of tilting objects in $\c H$ has a partial
order such that $T\leqslant T'$ if and only if $\ts{Fac}\, T\subseteq
\ts{Fac}\, T'$ where $\ts{Fac}\,T$ is the class of epimorphic images
of direct sums of copies of $T$. The Hasse diagram of this poset is denoted by
$\overrightarrow{\c K}_{\c H}$ and called the \emph{oriented graph of
  tilting objects in $\c H$} (see \cite{happel_unger3} for more details).

If $\c A$ is an additive category, then $\ts{ind}\, \c A$ denotes the
full subcategory of all indecomposable objects of $\c A$.

\paragraph{Galois coverings of $k$-categories}$\ $

Let $F\colon\c E\to \c B$ be a Galois covering with group $G$ between
$k$-categories (see the introduction). It is called \emph{connected}
if both $\c C$ and $\c B$ are connected and locally bounded.
Let $A$ be a connected and locally bounded $k$-category
and  $x_0\in \ts{ob}(A)$. A \textit{pointed Galois covering} $F\colon(\c
C,x)\to (A,x_0)$ is a connected Galois covering $F\colon\c C\to A$
endowed with  $x\in \ts{ob}(\c C)$ such that $F(x)=x_0$. A \textit{morphism
of pointed Galois coverings} $F\xrightarrow{u}F'$ from $F\colon(\c
C,x)\to (A,x_0)$ to $F'\colon(\c
C',x')\to (A,x_0)$ is a functor $u\colon \c C\to \c C'$ such that
$F'\circ u=F$ and $u(x)=x'$. Note that, given $F$ and $F'$, there is at
most one such morphism (see \cite[Lem. 3.1]{lemeur2}). This defines the
category $\ts{Gal}(A,x_0)$ of pointed Galois coverings. If $F\in
\ts{Gal}(A,x_0)$, then we let $F^{\to}$ be the full subcategory of
$\ts{Gal}(A,x_0)$ with objects those $F'$ such that there exists a morphism
$F\to F'$.

\paragraph{Covering properties on module categories
  (see \cite{bongartz_gabriel,riedtmann})}$\ $

 Let $F\colon \c E\to\c
B$ be a Galois covering with group $G$. The $G$-action on $\c E$
defines a $G$-action on $\ts{MOD}\,\c E$: If $M\in \ts{MOD}\,\c E$ and
$g\in
G$, then $^gM:=F\circ g^{-1}\in 
\ts{MOD}\,\c E$. If $X\in \ts{MOD}\,\c E$, \textit{the stabiliser of $X$} is the
subgroup $G_X:=\{g\in G\ |\ ^gX\simeq X\}$ of $G$. The Galois covering
$F$ defines two
exact functors: The extension-of-scalars functor $F_{\lambda}\colon
\ts{MOD}\,\c E\to \ts{MOD}\,\c B$ which is called the
\textit{push-down functor} and the restriction-of-scalars functor
$F.\colon \ts{MOD}\,\c B\to \ts{MOD}\,\c E$ which is called the \textit{pull-up
functor}. They form an adjoint pair $(\l,F_.)$ and $\l$ is
$G$-invariant, that is, $\l\circ g=\l$ for every $g\in G$. We refer the
reader to
\cite{bongartz_gabriel} for details one $\l$ and $F_.$.
For any $M,N\in \ts{mod}\,\c E$, the following maps
  induced by
  $F_{\lambda}$ are bijective:
 \begin{equation}
 \bigoplus\limits_{g\in G}\ts{Hom}_{\c E}(\,^gM,N)\to \ts{Hom}_{\c
   B}(F_{\lambda}M,F_{\lambda}N)\ \text{and}\ \bigoplus\limits_{g\in G}\ts{Hom}_{\c E}(M,\,^gN)\to \ts{Hom}_{\c
   B}(F_{\lambda}M,F_{\lambda}N)\ .\notag
\end{equation}
 An indecomposable
module  $X\in \ts{mod}\,\c B$
is called \textit{of the first kind
with respect to $F$} if and only if $\l\widetilde{X}\simeq X$ for some $\widetilde{X}\in \ts{mod}\,\c E$
(necessarily indecomposable). In such a case, one may choose
$\widetilde{X}$ such that $\l \widetilde{X}=X$. More
generally, $X\in \ts{mod}\,\c B$ is called of the first kind with respect to $F$ if
and only if it is the direct sum of indecomposable $\c B$-modules of
the first kind with respect to $F$.

\section{Covering techniques on the bounded derived category}
\label{sec:s2}$\ $

Let $F\colon\c C\to A$ be a Galois covering with group $G$ and with
$\c C$ and $A$ locally bounded categories of finite global
dimension. The $G$-action on $\ts{mod}\, \c C$ naturally defines
a $G$-action on $\dbc$, still denoted by $(g,M)\mapsto\, ^gM$, by triangle automorphisms. We 
introduce an exact functor
$\l\colon\dbc\to\dba$ induced by
$\l\colon\ts{mod}\, \c C\to\ts{mod}\, A$.
\begin{prop}
\label{prop1.1}
  There exists an exact functor
  $\l\colon\dbc\to\dba$ such that the following diagram commutes:
\begin{equation}
\label{eq:D1}
\xymatrix{
\ts{mod}\, \c C  \ar@{^(->}[r] \ar@{->}[d]_{\l} & \dbc
\ar@{->}[d]^{\l}\\
\ts{mod}\, A  \ar@{^(->}[r]  & \dba&.}\notag
\end{equation}
The functor $\l\colon\dbc\to\dba$ has the covering property, that is,
it is $G$-invariant and the two following
maps are linear bijections for every $M,N\in\dbc$:
\begin{align}
  &\bigoplus\limits_{g\in G}\dbc(\,^gM,N)\xrightarrow{\l} \dba(\l
  M,\l N) \ ,\notag\\
and\ & \bigoplus\limits_{g\in G}\dbc(M,\,^gN)\xrightarrow{\l} \dba(\l
  M,\l N)\ .\notag
\end{align}
\end{prop}
\noindent{\textbf{Proof:}} The existence and exactness of $\l\colon \dbc\to\dba$
follow from the exactness of $\l\colon\ts{mod}\, \c C\to \ts{mod}\, A$. On the other hand, 
 $\l$ induces an additive functor
$\l\colon \c K^b(\ts{mod}\, \c C)\to\c K^b(\ts{mod}\, A)$ between bounded homotopy
categories of complexes. It easily checked that it has the covering
property in the sense of the proposition. Since $A$ and $\c C$ have
finite global dimension, we deduce that $\l\colon\dbc\to\dba$ has the
covering property.\sq

\begin{rem}
  \label{rem1.2}
It follows from the preceding proposition that $\l\colon \dbc\to \dba$
is faithful.
\end{rem}

We are mainly interested in indecomposable objects $X\in\dba$ which
are of the form $\l\w X$ for some $\w X\in\dbc$. The following shows
that the possible objects $\w X$ lie in the same $G$-orbit for a given $X$.
\begin{lem}
\label{lem1.3}
  Let $X,Y\in\dbc$ be such that $\l X$ and $\l Y$ are indecomposable
  and isomorphic in
  $\dba$. Then $X\simeq\,^g Y$ for some
  $g\in G$.
\end{lem}
\noindent{\textbf{Proof:}} Let $u\colon \l X\to \l Y$ be an
isomorphism in $\dba$. By \ref{prop1.1}, there exists
 $(u_g)_g\in
\bigoplus\limits_{g\in G}\dbc(X,\,^gY)$ such that $u=\sum\limits_{g\in
  G}\l(u_g)$. Since $\l X$ and $\l
Y$ are indecomposable, there exists $g_0\in G$ such that $\l
(u_{g_0})\colon\l X\to \l Y$ is an isomorphism. Since $\l\colon
\dbc\to \dba$ is exact and faithful, $u_{g_0}\colon
X\to \,^{g_0}Y$ is an isomorphism in $\dbc$.\sq

\section{Transforming tilting complexes into tilting modules}
\label{sec:s3}$\ $

Let $\c H$ be a hereditary abelian category over $k$ with finite
dimensional $\ts{Hom}$-spaces, split idempotents and tilting
objects. Let $n$  the rank of its
Grothendieck group. For short, we set $\ts{Hom}:=\ts{Hom}_{\dbh}$ and
$\ts{Ext}^i(X,Y):=\ts{Hom}_{\dbh}(X,Y[i])$.
We write $\c T$ for the class of
objects $T\in\dbh$ such that:
\begin{enumerate}
\item[(a)] $T$ is multiplicity-free and has $n$ indecomposable summands.
\item[(b)] $\ts{Ext}^i(T,T)=0$ for every $i\geqslant 1$.
\end{enumerate}
We identify an object in $\c T$ with its isomorphism class.
A complex $T$ lies in $\c T$ if and only if $T[1]\in\c T$. Also, all tilting
complexes in $\dbh$ and, therefore, all tilting objects in $\c H$, lie
on $\c T$.
Given $T\in\dbh$, we let $\left<T\right>$ be the
smallest full subcategory of $\dbh$ containing $T$ and stable under
direct sums, direct summands and shifts in either direction.
The aim of this section is to define elementary transformations on
objects in $\c T$ which, by repetition, allow one to relate any two
objects in $\c T$.
For this purpose, we introduce some notation. Given $T\in \c T$, we
have a unique decomposition $T=Z_0[i_0]\oplus
Z_1[i_0+1]\oplus\ldots\oplus Z_l[i_0+l]$ where
$Z_0,\ldots,Z_l\in \c H$ and $Z_0,Z_l\neq 0$. Here, $Z_i$ needs not
be indecomposable. We
let $r(T)$ be the number of indecomposable summands of
$Z_1[i_0+1]\oplus\ldots\oplus Z_l[i_0+l]$. Note that:
 $r(T)\in\{0,\ldots,n-1\}$;
 $r(T)=0$ if and only if $T[-i_0]$ is a tilting object in $\c
  H$; and
 $r(T)=r(T[1])$. We are interested in transformations which map an
 object $T\in\c T$ to $T'$ such that $r(T')<r(T)$. Hence, by repeating
 the process, we may end up with a tilting object in $\c H$ (up to a
 shift). 

\paragraph{Transformations of the first kind}$\ $

Our first elementary transformation is given by the following lemma.
\begin{lem}
\label{lem2.1}
  Let $T\in\c T$. There exists $T'\in\c T$ such that
$T'\in \left<T\right>$,
  $r(T')\leqslant r(T)$
 and  $T'=Z_0'\oplus Z_1'[1]\oplus\ldots\oplus Z_{l'}'[l']$
    where:
    \begin{enumerate}
    \item[(a)] $Z_0',\ldots,Z_{l'}'\in \c H$ and $Z_0',Z_{l'}'\neq 0$.
    \item[(b)] $\ts{Hom}(Z_0',Z_1'[1])\neq 0$ if $l'\neq 0$.
    \end{enumerate}
\end{lem}
\noindent{\textbf{Proof:}} Given $T'\in \left<T\right>$, we have the
unique decomposition $T'=Z_0'[i_0']\oplus
Z_1'[i_0'+1]\oplus\ldots\oplus Z_l'[i_0'+l']$ as explained at the
beginning of the section. We choose $T'\in \left<T\right>\cap \c T$ such that
 $r(T')\leqslant r(T)$ and such that the pair $(l',r(T'))$ is minimal for the
 lexicographical order. We may assume that $i_0=0$. We prove that $T'$
 satisfies (a) and (b). If $l'=0$, there is nothing to prove. So we
assume that $l'>0$. Assume first that $Z_1'=0$. Then we let $T''$ be
as follows:
\begin{equation}
  T'':=Z_0'\oplus Z_2'[1]\oplus Z_3'[2]\oplus\ldots \oplus Z_{l'}'[l'-1]\ .\notag
\end{equation}
Then  $T''\in \left<T'\right>=\left<T\right>$. Also,
$\ts{Ext}^i(T'',T'')=0$ for every $i\geqslant 1$ because 
$T'\in\c T$ and $\c H$ is hereditary. Finally,
 $T''$ is the direct sum of $n$ pairwise non isomorphic
  indecomposable objects. Thus, $T''\in\left<T\right>\cap\c T$ and
  $(l'-1,r(T''))<(l',r(T'))$ which contradicts the minimality of
  $(l',r(T'))$.  So $Z_1'\neq 0$. Now, assume
that $\ts{Hom}(Z_0',Z_1'[1])=0$. We let $T''$ be the
following object:
\begin{equation}  
T'':=(Z_0'\oplus Z_1')\oplus Z_2'[2]\oplus Z_3'[3]\oplus\ldots\oplus Z_{l'}'[l']\ .\notag
\end{equation}
As above, we have $T''\in \left<T\right>\cap \c
T$ and $(l',r(T''))<(l',r(T'))$ which contradicts the minimality of
$(l',r(T'))$. So $\ts{Hom}(Z_0',Z_1'[1])\neq
0$.\sq

With the notations of \ref{lem2.1}, we say that $T$ and
$T'$ are related by a \textit{transformation of the first kind}.

\paragraph{Transformations of the second kind}$\ $

We now turn to the second elementary transformation. It is inspired
by the characterisation of the quiver of tilting objects in hereditary
categories (see
\cite{happel_unger3} and also \cite{bmrrt} for the
corresponding construction in cluster categories). 
Let $T,T'\in\c T$
be such that $T=X\oplus\o T$ with $X$
indecomposable, $T'=Y\oplus\o T$ with $Y$ indecomposable and there
exists a triangle $X\xrightarrow{u}M\xrightarrow{v}Y\to X[1]$ such
that $u$ is a left minimal $\ts{add}(\o T)$-approximation or $v$ is
a right minimal $\ts{add}(\o T)$-approximation. In such a situation, we
 say that $T$ and $T'$ are related by a \textit{transformation of
the second kind}.
\begin{rem}
  \label{rem2}
Following \cite{happel_unger3}, if $T\to T'$ is an arrow in
$\overrightarrow{\c K}_{\c H}$ then $T$ and $T'$ are relatd by a
transformation of the second kind.
\end{rem}

 Note that, with the previous notations, both $u$ and $v$ are minimal $\ts{add}(\o
T)$-approximations, as shows the following lemma.
\begin{lem}
\label{lem2.7}
  Let $T\in \c T$. Assume that $T=X\oplus \overline{T}$ with
  $X$ indecomposable.
\begin{enumerate}
\item[(a)]
Let  $X\xrightarrow{u}M\xrightarrow{v}Y\to X[1]$ be a triangle where $u$ is
  a left minimal $\ts{add}(\o T)$-approximation. Then $v$ is a right minimal
   $\ts{add}(\o T)$-approximation.
\item[(b)]
Let  $Y\xrightarrow{u}M\xrightarrow{v}X\to Y[1]$ be a triangle where $v$ is
  a right minimal $\ts{add}(\o T)$-approximation. Then $u$ is a left minimal
   $\ts{add}(\o T)$-approximation.
\end{enumerate}
\end{lem}
\noindent{\textbf{Proof:}} We only prove (a) because the
proof of (b) is similar. Every morphism $\o T\to Y$ factorises through
$v$ because $\ts{Hom}(\o T,X[1])=0$. So $v$ is a right $\ts{add}(\o
T)$-approximation. Let $\alpha\colon M\to M$ be a morphism such that
$v\alpha=v$. So there exists $\lambda\colon M\to X$ such that
$u\lambda=\alpha-\ts{Id}_M$. Note that $u$ is not a section because $T$ is
multiplicity-free. So  $u \lambda$ is
nilpotent and $\alpha=\ts{Id}_M+u \lambda$ is an
isomorphism. Therefore $v$ is right minimal.\sq

It is not true that any two objects $T,T'\in\c T$ can be
related by a sequence of transformations of second kind (whereas this
is the case, for example, for tilting objects in a cluster category,
see \cite{bmrrt}). However,
we have the following result from \cite{happel_unger3}.
\begin{prop}
\label{prop2.7}
  Assume that  at least one of the two following assertions is true:
\begin{enumerate}
\item[(a)] $\c H=\ts{mod}\,kQ$ where $Q$ is a finite
  connected quiver without oriented cycles and of Dynkin type.
\item[(b)] $\c H$ has no non-zero projective object and $\c D^b(\c H)$
  is triangle equivalent to $\c D^b(\ts{mod}\, kQ)$ with $Q$ a
  connected finite quiver without oriented cycles.
\end{enumerate}
Then $\overrightarrow{\c K}_{\c H}$ is connected. In partiular (see
\ref{rem2}) for every tilting objects $T,T'\in\c H$. There exists a
sequence $T=T_0,\ldots,T_l=T'$ of
  tilting objects in $\c H$ such that $T_i$ and $T_{i+1}$ are related
  by a transformation of the second kind for every $i$. 
\end{prop}

We are going to prove that any $T\in\c T$ can be related to some
tilting object in $\c H$ by a sequence of transformations of the first
or of the second kind.
Let $T\in\c T$. With the notations established at the beginning of the
section, assume that
$\ts{Hom}(Z_0,Z_1[1])\neq 0$. Since the ordinary quiver of
$\ts{End}(T)$ has no oriented cycle, there exists $M\in \ts{add}(Z_1[i_0+1])$ 
indecomposable such that:
\begin{enumerate}
\item $\ts{Hom}(Z_0[i_0],M)\neq 0$.
\item $\ts{Hom}(Z,M)=0$ for any indecomposable direct summand $Z$ of
  $\bigoplus\limits_{t=1}^lZ_t[i_0+t]$ not isomorphic to $M$.
\end{enumerate}
 Let
$\o T$ be such that $T=\o T\oplus M$. Let $B\to M$ be a right minimal
$\ts{add}(\o T)$-approximation of $M$. Complete it into a triangle in
$\dbh$:
\begin{equation}
  M^*\to B\to M\to M^*[1]\ .\tag{$\Delta$}
\end{equation}
\begin{lem}
\label{lem2.2}
With the above setting, let $T'=\o T\oplus M^*$. Then $T'\in \c T$
and $T,T'$ are related by a transformation of the second kind. Moreover, $r(T')<r(T)$. 
\end{lem}
\noindent{\textbf{Proof:}} We only need to prove that $T'\in\c T$. We may assume
that $i_0=0$. By hypothesis on $M$, we have $B\in \ts{add}(Z_0)\subseteq
\c H$. Since $M\in\c H[1]$, the triangle $\Delta$ does not split. We
now list some properties on $T'$. In most cases, the proof is due to
arguments taken from \cite[§ 6]{bmrrt}. Although these arguments
were originally given in the setting of cluster categories (that is,
triangulated categories which are Calabi--Yau of dimension $2$), it is
easily verified that they still work in our situation (that is, the
Calabi--Yau property is unnecessary):
\begin{enumerate}
\item $\ts{Ext}^1(\o T,M^*)=0$ (see \cite[Lem. 6.3]{bmrrt}).
\item $\ts{Ext}^i(\o T,M^*)=0$ for every $i\geqslant 2$. Indeed,
  applying $\ts{Hom}(\o T,-)$  to $\Delta$ gives the
  exact sequence
  $$0=\ts{Ext}^{i-1}(\o T,M)\to \ts{Ext}^i(\o T,M^*)\to \ts{Ext}^i(\o T,B)=0\ .\notag$$
\item $\ts{Ext}^i(M^*,\o T)=0$ for every $i\geqslant 1$. Indeed, applying $\ts{Hom}(-,\o T)$ to
  $\Delta$ gives the exact sequence $$0=\ts{Ext}^i(B,\o T)\to \ts{Ext}^i(M^*,\o
  T)\to \ts{Ext}^{i+1}(M,\o T)=0\ .\notag$$
\item The map $M^*\to B$ is a left minimal $\ts{add}(\o
  T)$-approximation (see \cite[Lem. 6.4]{bmrrt}).
\item $M^*$ is indecomposable and does not lie on $\ts{add}(\o T)$ (see
  \cite[Lems. 6.5, 6.6]{bmrrt}). Therefore $T'$ is the direct sum of
  $n$ pairwise indecomposable objects.
\item $M^*\in\c H$. Indeed, $M$ is indecomposable and there are two non-zero maps $M[-1]\to
  M^*$ and $M^*\to B$ with $M[-1],B\in\c H$.
\item $\ts{Ext}^1(M^*,M^*)=0$ (see
  \cite[Lem. 6.7]{bmrrt}).
\item $\ts{Ext}^i(M^*,M^*)=0$ for every $i\geqslant 2$ because $M^*$ is
  indecomposable and $\c H$ is hereditary.
\end{enumerate}
The facts $1.-8.$ express that $T'\in\c T'$. Moreover, $r(T')<r(T)$
because $M^*\in\c H$ and $M\in\c H[1]$.\sq

\begin{lem}
  \label{lem2.3}
Let $T\in\c T$. Let $\c A$ be the smallest subclass of $\c T$ containing
$T$ and stable under transformations of
the first or of the second kind. Then $\c A$
contains a tilting object in $\c H$.
\end{lem}
\noindent{\textbf{Proof:}} Let $T'\in\c A$ be such that $r(T')$ is minimal
for this property. Assume that $r(T')>0$. By
\ref{lem2.1} and \ref{lem2.2}, there exists
 $T''\in \c A$ such that $r(T'')<r(T')$. This
contradicts the minimality of $r(T')$. Hence $r(T')=0$ and there exists
an integer $i_0$ such that $T'[-i_0]$ is a tilting object in
$\c H$ and lies in $\c A$.\sq

 The following result is a direct
consequence of \ref{prop2.7} and \ref{lem2.3}.
\begin{prop}
  \label{prop2.4}
Assume that  $\overrightarrow{\c K}_{\c H}$ is connected. Let $T\in\c T$. Then
$\c T$ is the smallest subset of $\c T$ containing $T$ and stable
under transformations of the first or the second kind.
\end{prop}

\begin{rem}
  \label{rem2.5}
\begin{enumerate}
\item[(a)] A tilting object in $\c H$ generates $\dbh$. By definition
of the two kinds of transformations, 
\ref{prop2.4} implies, under the same hypotheses, that any
$T\in\c T$ generates $\dbh$.
\item[(b)] Let $A$ be an algebra derived equivalent to a hereditary
  algebra. Then \ref{prop2.7} implies that the
  conclusion of 
  \ref{prop2.4} holds true if one replaces $\c H$ by $\ts{mod}\,
  A$.
\end{enumerate}
\end{rem}

\section{Tilting complexes of the first kind}
\label{sec:s4}$\ $

Throughout this section, we assume that $A$ is an algebra derived
equivalent to a hereditary abelian category $\c H$ such that
$\overrightarrow{\c K}_{\c H}$  is connected.
 We denote by $n$ the rank of its
Grothendieck group and $\Theta\colon \c D^b(\c H)\to \dba$ a triangle equivalence. We fix
a Galois covering $F\colon \c C\to A$ with group $G$ and with $\c
C$ locally bounded. We use \ref{prop1.1} without
reference.
The aim of this section is to prove that the
following facts hold true for any tilting complex $T\in\dba$:
\begin{enumerate}
\item[($H_1$)] For every indecomposable direct summand $X$ of $T$,
  there exists $\w X\in\dbc$ such that $\l\w X\simeq X$ in $\dba$.
\item[($H_2$)] $\w X\not\simeq\,^g\w X$ for every indecomposable direct
  summand $X$ of $T$ and $g\in G\backslash\{1\}$.
\item[($H_3$)] If $\psi\colon
A\xrightarrow{\sim} A$ is an automorphism such that $\psi(x)=x$ for
every $x\in \ts{ob}(A)$, then $\psi_{\lambda}X\simeq X$ in $\dba$ for
every indecomposable direct summand $X$ of $T$.
\end{enumerate}
Some results presented in this section have been proved in \cite[§
3]{lemeur5} in the case of tilting modules.

\paragraph{Proof of assertion ($H_1$)}$\ $

In this paragraph, we  prove the following.
\begin{prop}
\label{prop3.1}
Let $A$ be as above. Let $T\in \dba$ be a tilting complex. Then:
\begin{enumerate}
\item[(a)] For every indecomposable direct
summand $X$ of $T$ there exists $\w X\in\dbc$ (necessarily
indecomposable) such that $\l\w X\simeq X$.
\end{enumerate} Moreover, the class
$\{\w X\in\dbc\ |\ \text{$\l\w X$ is an indecomposable direct summand
    of $T$}\}$ satisfies the following:
  \begin{enumerate}
  \item[(b)] It generates the triangulated category $\dbc$.
  \item[(c)] It is stable under the action of $G$.
  \item[(d)] $\dbc(\w X,\,^g\w Y[i])=0$ for every $\w X,\w Y$ in this class, $i\neq 0$
    and $g\in G$.
  \end{enumerate}
\end{prop}

We need the
two following dual lemmas.
  \begin{lem}
\label{lem3.2}
Let $\Delta:X\xrightarrow{u} M\to Y\to X[1]$
be triangle in $\dba$ such that:
\begin{enumerate}
\item[(a)] There exists $\w X\in\dbc$ satisfying $X=F_{\lambda}\w{X}$.
\item[(b)] $M=M_1\oplus \ldots \oplus M_t$ where $M_1,\ldots,M_t$
  are indecomposables such that there exist indecomposable objects $\w
  M_1,\ldots,\w M_t$ satisfying $\l\w M_i=M_i$ for every $i$.
\item[(c)] $\dba(Y,M[1])=0$.
\end{enumerate}
Then $\Delta$ is isomorphic to a triangle in $\dba$:
\begin{equation}
X \xrightarrow{\begin{bmatrix} F_{\lambda}u_1'& \ldots &
    F_{\lambda}u_t'\end{bmatrix}^t} M_1\oplus\ldots\oplus
    M_t\to Y\to X[1]\notag
\end{equation}
where $u'_i\in {\dbc}(\w{X},\,^{g_i}\w{M}_i)$ for some
$g_i\in G$ for every $i$.
\end{lem}

\begin{lem}
\label{lem3.3}
Let $\Delta:X\to M\xrightarrow{v} Y\to X[1]$
be triangle in $\dba$ such that:
\begin{enumerate}
\item[(a)] There exists $\w Y\in\dbc$ satisfying $Y=F_{\lambda}\w{Y}$.
\item[(b)] $M=M_1\oplus \ldots \oplus M_t$ where $M_1,\ldots,M_t$
  are indecomposables such that there exist indecomposable objects $\w
  M_1,\ldots,\w M_t$ satisfying $\l\w M_i=M_i$ for every $i$.
\item[(c)] $\dba(M,X[1])=0$.
\end{enumerate}
Then $\Delta$ is isomorphic to a triangle in $\dba$:
\begin{equation}
X\to  M_1\oplus\ldots\oplus
    M_t\xrightarrow{\begin{bmatrix} F_{\lambda}v_1'& \ldots &
    F_{\lambda}v_t'\end{bmatrix}} Y\to X[1]\notag
\end{equation}
where $v'_i\in {\dbc}(\,^{g_i}\w{M}_i,\w{Y})$ for some
$g_i\in G$ for every $i$.
\end{lem}

\noindent{\textbf{Proof of \ref{lem3.2}}}: 
We say that a morphism $u\in \dba(X,M_i)$ is homogeneous of
degree $g\in G$ if and only if there exists $u'\in
\dbc(\w{X},\,^g\w{M}_i)$ such that $u=\l(u')$. Since
$\l\colon\dbc\to \dba$ has the covering property, any morphism $X\to
M_i$ is (uniquely) the
sum of $d$ non zero homogeneous morphisms of pairwise different degrees (with
$d\geqslant 0$).
Let
    $u=\begin{bmatrix} u_1& \ldots &
    u_t\end{bmatrix}^t$ with $u_i\colon X\to M_i$ for each $i$. We may
  assume that $u_1\colon X\to M_1$ is not
    homogeneous. Thus $u_1=h_1+\ldots+h_d$,
where $d\geqslant 2$ and $h_1,\ldots,h_d\colon X\to M_1$ are non-zero
homogeneous morphisms of pairwise different degrees. In order to prove
the lemma, it
suffices to prove that $\Delta$ is isomorphic to a triangle
$X\xrightarrow{u'}M\to Y\to X[1]$ with
$u'=\begin{bmatrix}u_1'&u_2&\ldots&u_t\end{bmatrix}^t$ such that
 $u_1'$ is equal to the sum of at most $d-1$ non-zero homogeneous morphisms
$X\to M_1$ of pairwise different degrees.
For simplicity we adopt the following notations:
\begin{enumerate}
\item $\overline{M}=M_2\oplus\cdots\oplus M_t$ (so
  $M=M_1\oplus \overline{M}$).
\item $\overline{u}=\begin{bmatrix} u_2& \ldots &
    u_t\end{bmatrix}^t\colon X\to \overline{M}$ (so $u=\begin{bmatrix}
    u_1&\overline{u} \end{bmatrix}^t\colon X\to M_1\oplus\overline{M}$).
\item $\overline{h}=h_2+\ldots+h_d\colon X\to M_1$ (so
  $u_1=h_1+\overline{h}$).
\end{enumerate}
Applying the functor ${\dba}(-,M_1)$  to $\Delta$ gives the exact sequence:
\begin{equation}
\dba(M_1\oplus\overline{M},M_1)\xrightarrow{\ts{Hom}(u,M_1)} \dba(X,M_1) \to
\dba(Y,M_1[1])=0\ .\notag
\end{equation}
So there exists $[\lambda,\mu]\colon M_1\oplus\overline{M}\to M_1$ such
that $h_1=[\lambda,\mu]u$.
Hence:
\begin{equation}
h_1=\lambda u_1+\mu\overline{u}=\lambda
h_1+\lambda\overline{h}+\mu\overline{u}\ .\tag{$i$}
\end{equation}
We distinguish two cases according to whether $\lambda\in \ts{End}_{\dba}(M_1,M_1)$ is
invertible or nilpotent.
 If $\lambda$ is invertible, then the following is an
isomorphism in $\dba$:
\begin{equation}
\theta:=\begin{bmatrix}
\lambda & \mu\\
0 & \ts{Id}_{\overline{M}}
\end{bmatrix}\colon M_1\oplus \overline{M}
\to M_1\oplus \overline{M}\ .\notag
\end{equation}
Using $(i)$ we deduce an isomorphism of triangles:
\begin{equation}
\xymatrix{
 X \ar@{->}[r]^{\begin{bmatrix} u_1&
    \overline{u}\end{bmatrix}^t} \ar@{=}[d] & M_1\oplus \overline{M}
\ar@{->}[r] \ar@{->}[d]^{\theta} & Y \ar@{->}[r]
\ar@{->}[d]^{\sim}& X[1]  &\Delta\\
 X \ar@{->}[r]^{\begin{bmatrix} h_1&
    \overline{u}\end{bmatrix}^t} & M_1\oplus\overline{M}\ar@{->}[r] &
Y\ar@{->}[r] & X[1] & \Delta'\ .
}\notag
\end{equation}
Since $h_1\colon X\to M_1$ is homogeneous, $\Delta'$ satisfies the
our
requirements.
 If $\lambda$ is nilpotent, let $p\geqslant 0$
be such that $\lambda^p=0$. Using $(i)$ we get the following
equalities:
\begin{equation}
\begin{array}{rclc}
h_1&=&\lambda^2h_1+(\lambda^2+\lambda)\overline{h}+(\lambda+\ts{Id}_{M_1})\mu\overline{u}\\
\vdots&\vdots&\vdots\\
h_1&=&\lambda^th_1+(\lambda^t+\lambda^{t-1}+\ldots+\lambda)\overline{h}+(\lambda^{t-1}+\ldots+\lambda+\ts{Id}_{M_1})\mu\overline{u}\\
\vdots&\vdots&\vdots\\
h_1&=&\lambda^ph_1+(\lambda^p+\lambda^{p-1}+\ldots+\lambda)\overline{h}+(\lambda^{p-1}+\ldots+\lambda+\ts{Id}_{M_1})\mu\overline{u}&.
\end{array}\notag
\end{equation}
Since $\lambda^p=0$ and $u_1=h_1+\overline{h}$ we infer that:
\begin{equation}
u_1=\lambda'\overline{h}+\lambda'\mu\overline{u}\ ,\notag
\end{equation}
where $\lambda'$ is the invertible morphism
$\ts{Id}_{M_1}+\lambda+\ldots+\lambda^{p-1}\in \ts{End}_{\dba}(M_1,M_1)$. So
we have an isomorphism in $\dba$: 
\begin{equation}
\theta:=\begin{bmatrix}
\lambda'&\lambda'\mu\\
0&\ts{Id}_{\overline{M}}
\end{bmatrix}\colon M_1\oplus\overline{M}\to
M_1\oplus\overline{M}\ .\notag
\end{equation}
Consequently we have an isomorphism of triangles:
\begin{equation}
\xymatrix{
 X \ar@{=}[d]^{}
 \ar@{->}[r]^{\begin{bmatrix}
\overline{h}&
\overline{u}
\end{bmatrix}^t} &M_1\oplus\overline{M} \ar@{->}[d]^{\theta}
\ar@{->}[r]& Y\ar@{->}[r] \ar@{->}[d]^{\sim}& X[1] & \Delta'\\
 X\ar@{->}[r]^{
\begin{bmatrix}
u_1&
\overline{u}
\end{bmatrix}^t
} & M_1\oplus \overline{M} \ar@{->}[r] & Y \ar@{->}[r] & X[1]&
\Delta
}\notag
\end{equation}
where $\overline{h}=h_2+\ldots+h_p$ is the sum of $p-1$ non zero homogeneous
morphisms of pairwise different degrees.
So $\Delta'$ satisfies our requirements.\sq

The proof of \ref{lem3.3} is the dual of the one of
\ref{lem3.2} so we omit it.
Now we
can prove \ref{prop3.1}.\\
\noindent{\textbf{Proof of \ref{prop3.1}}:} 
If (a) holds true, then so does (c)  because $\l\colon \dbc\to \dba$
is $G$-invariant.
Recall that $\Theta\colon\dbh\to\dba$ is a triangle equivalence. As in
Section~\ref{sec:s3}, we write $\c T$ (or $\c T'$) for the set of isomorphism
classes of objects $T\in\dbh$ (or $T\in\dba$) such that
 $T$ is the direct sum of $n$ pairwise non isomorphic
  indecomposable objects and
 $\dbh(T,T[i])=0$ (or $\dba(T,T[i])=0$, respectively) for every $i\geqslant 1$.
Therefore:
\begin{enumerate}
\item[($i$)] $\Theta$ defines a bijection $\Theta\colon\c T\to\c
T'$. Under this bijection, tilting complexes in $\dbh$ correspond to
tilting complexes in $\dba$. 
\end{enumerate}
 We  prove that (a) and (b) hold true for any  $T\in\c T'$ (and therefore for any
tilting object in $\dba$). For this purpose, we use the results of
Section~\ref{sec:s3}. First of all, remark that:
\begin{enumerate}
\item[($ii$)] The assertions (a) and (b) hold true for $T=A$. In this case, $\l\w X$ is an
indecomposable summand of $A$ if and only if $\w X$ is an
indecomposable projective $\c C$-module.
\end{enumerate}
If $\w X\in\dbh$, then
$\l(\w X[1])=(\l\w X)[1]$. Therefore:
\begin{enumerate}
\item[($iii$)] Let $T,T'\in\c T'$ be such that $\Theta^{-1}(T')$ is obtained from $\Theta^{-1}(T)$
  by a transformation of the first kind. Then (a) and (b) hold
  true for $T$ if and only if they do so for $T'$.
\end{enumerate}
Now assume that $T,T'\in\c T'$ are such that $\Theta^{-1}(T')$ is
obtained from $\Theta^{-1}(T)$ by a transformation of the second
 kind. We prove that (a) and
(b)  hold true for $T$ if and
only if they do so for $T'$. In such a situation there exist
$X,Y\in\dba$ indecomposables and $\overline{T}\in\dba$ such
that $T=X\oplus\overline{T}$ and $T'=Y\oplus\overline{T}$. Also, there
exists a triangle in $\dba$ of one the two following forms:
\begin{enumerate}
\item $X\to M\to
  Y\to X[1]$ where $M\in\ts{add}(\overline{T})$.
\item $Y\to M\to
  X\to Y[1]$ where $M\in\ts{add}(\overline{T})$.
\end{enumerate}
Assume that (a) and (b) hold true for $T$ and that there is a triangle
$X\to M\to Y\to X[1]$ (the other cases are dealt with using similar
arguments). In order to prove that (a) and (b) hold
true  for $T'$ 
we prove that $Y\simeq \l\w Y$ for some $\w Y\in\dbc$. Fix an indecomposable decomposition
$M=M_1\oplus\ldots\oplus M_t$. By assumption on $T$,
there exist indecomposable objects 
$\w X,\w M_1,\ldots,\w M_t\in\dbc$ such that $\l\w X\simeq X,\l\w
M_1\simeq M_1,\ldots,\l\w M_t\simeq M_t$. Using these isomorphisms, we
identify $\l\w X$ and $\l\w M_i$ to $X$ and $M_i$, respectively. By \ref{lem3.2}, 
there exist $g_1,\ldots,g_t\in G$ and morphisms
 $u_i\in \dbc(\w{X},\,^{g_i}\w{M_i})$
 (for $i\in\{1,\ldots,t\}$) such that the triangle $X\to M\to Y\to
 X[1]$ is isomorphic to a triangle of the form:
 \begin{equation}
 X\xrightarrow{
   \begin{bmatrix}
     \l(u_1)& \ldots&  \l(u_t)
   \end{bmatrix}^t}
M\to Y\to X[1]\ .\notag
 \end{equation}
Set $u=
\begin{bmatrix}
  u_1& \ldots& u_t
\end{bmatrix}^t\colon \w X\to \w M_1\oplus\ldots\oplus \w M_t$. We
complete $u$ into a triangle $\w X\xrightarrow{u}\w
M_1\oplus\ldots\oplus \w
M_t\xrightarrow{v} \w X[1]$ in $\dbc$. So we have a triangle
$X\xrightarrow{\l(u)}M\xrightarrow{\l(v)}\l \w Y\to X[1]$ in $\dba$. Therefore
$Y\simeq \l\w Y$. So (a) holds for $T'$ and the class $\{\w Z\ |\
\text{$\l\w Z$ is an
  indecomposable  direct
summand of $T'$}\}$ coincides with the class  $\{\,^g\w Y\ |\ g\in G\}\cup\{\w Z\
|\ \text{$\l\w Z$ is an indecomposable direct summand of
  $\overline{T}$}\}$ (see \ref{lem1.3}). Because (b)
holds true for $T$ and
because of the triangle  $\w
X\to \w M_1\oplus\ldots\oplus\w M_t\to \w Y\to \w X[1]$, we deduce
that (b) holds for $T'$. So we have proved the following:
\begin{enumerate}
\item[($iv$)] Let $T,T'\in\c T'$ be such that $\Theta^{-1}(T')$ is
  obtained from $\Theta^{-1}(T)$ by a  transformation 
  of the second kind. Then (a) and (b) hold true for $T$ if and only if
  they do so for $T'$.
\end{enumerate}
By \ref{prop2.4} and ($i-iv$), the assertions (a), (b) and
(c) are
 satisfied for any $T\in\c T$.
Finally, if $T$ is a tilting complex, then (d) follows from
the fact that $\dba(T,T[i])=0$ for every $i\neq 0$ and from \ref{prop1.1}.
\sq

It is interesting to note that the
transformations of the second kind in $\dba$ give rise to
transformations of the second kind in $\dbc$. Indeed, 
let $T,T'$ be in $\c T'$ where $\c T'$ is as in the proof of
\ref{prop3.1}. Assume that $T=M\oplus\o T$ with $M$
indecomposable, $T=M^*\oplus\o T$ with $M^*$ indecomposable and there
exists a triangle $\Delta:\ M\to B\to
M^*\to M[1]$ in $\dba$ where $M\to B$  (or $B\to M^*$) is a
left minimal $\ts{add}(\o T)$-approximation of $M$ (or a right minimal
$\ts{add}(\o T)$-approximation of $M^*$, respectively). Then the following holds.
\begin{lem}
  \label{lem3.10}
Keep the above setting. Let $B=\bigoplus\limits_{i=1}^t\, B_i$
be an indecomposable decomposition (maybe with multiplicities). Then there exists a triangle
$\w{\Delta}:\ \w M\xrightarrow{u} \bigoplus\limits_{i=1}^t\,^{g_i}\w B_i\xrightarrow{v}
\,^{g_0}\w M^*\to M[1]$ in $\dbc$ whose image under $\l$ is isomorphic
to $\Delta$. Moreover, if  $\c X$ (or $\c X'$) denotes the additive
full subcategory of $\dbc$ generated by the indecomposables
$X\in\dbc$ not isomorphic to $\w M$ (or to $\w M^*$) and such that
$\l X$ is an indecomposable summand of $T$ (or of $T'$, respectively), then:
\begin{enumerate}
\item[(a)] $u$ is a left minimal
  $\c X$-approximation.
\item[(b)] $v$ is a right minimal
  $\c X'$-approximation.
\end{enumerate}
\end{lem}
\noindent{\textbf{Proof:}} The existence of $\w{\Delta}$ follows from
the proof of \ref{prop3.1}. So $\l(u)$ is a left minimal
$\ts{add}(\o T)$-approximation. This and the exactness of $\l$ imply that
$u$ is left minimal. Let $f\colon \w M\to \,^g\w Y$ be a non-zero
morphism where $^g\w Y\in\c X$ and $Y\in\ts{add}(T)$. The linear map $\bigoplus\limits_{h\in
  G}\dbc(\w M,\,^h\w M)\to\ts{End}_{\dba}(M,M)$ induced by $\l$ is bijective. Also
$\ts{dim}_k\ts{End}_{\dba}(M,M)=1$ because $M$ is an indecomposable and
$\dba(M,M[i])=0$ for every $i>0$. So $^g\w Y\not\simeq\,^h\w M$ for
every $h\neq 1$. Hence $Y\in\ts{add}(\o T)$ and, therefore, $\l(f)$ factorises
through $\l(u)$:
\begin{equation}
  \xymatrix{
M \ar[r]^{
  \l(u)} \ar[rd]_{\l(f)}& \bigoplus\limits_{i=1}^t B_i\ar[d]^{
f'} \\
& Y & .
}\notag
\end{equation}
There exists $(f'_h)_h\in\bigoplus\limits_{h\in
  G}\dbc(\bigoplus\limits_{i=1}^t\,^{g_i}\w B_i,\,^h\w Y)$ such that
$f'=\sum\limits_{h\in G}\l(f'_g)$ because of the covering property of
$\l$. So $\l(f)=\sum\limits_{h\in G}\l(f'_h u)$ and, therefore, $f=f'_gu$ for the
same reason. Thus $u$ is a left minimal $\c
X$-approximation. Similarly, $v$ is a right minimal $\c X'$-approximation.\sq

Since tilting $A$-modules are particular cases of tilting complexes,
we get the following result.
\begin{cor}
  \label{cor3.4}
Let $A$ be an algebra derived equivalent to a hereditary abelian
category $\c H$ such that $\overrightarrow{\c K}_{\c H}$ is
connected. Let
$F\colon \c C\to A$ be a Galois covering with group $G$ where $\c C$
is locally bounded, $T$ a tilting $A$-module and $X\in\ts{mod}\,
A$ an
indecomposable summand of $T$. Then there exists $\w X\in\ts{mod}\,\c C$
such that $\l\w X\simeq X$.
\end{cor}
\noindent{\textbf{Proof:}} 
By \ref{prop3.1}, such an $\w X$ exists in $\dbc$. We
prove that $\w X$ is isomorphic to a $\c C$-module. Let $P\in\ts{mod}\,\c C$
be 
projective and  $i\in\mathbb{Z}\backslash\{0\}$. Then $\l P\in\ts{mod}\,
A$ is projective and $\dba(\l P,X[i])=0$ because $X$ is an $A$-module.
On the other hand, the spaces $\dba(\l P,X[i])$ and
$\bigoplus\limits_{g\in G}\dbc(\,^gP,\w X[i])$ are isomorphic. So
$\dbc(P,\w X[i])=0$ for every $i\neq 0$. Thus,
$\w X\simeq H^0(\w X)\in\ts{mod}\,\c C$.\sq

\paragraph{Proof of assertion ($H_2$)}$\ $

\begin{prop}
  \label{prop3.5}
Let $A$ be as in \ref{cor3.4},
$F\colon \c C\to A$ a Galois covering with group $G$ where $\c C$
is locally bounded and $X\in\ts{ind}\, A$ a direct
summand of a tilting complex in $\dba$. 
Assume that $\l\w X\simeq X$ for some $\w X\in\dbc$. Then $^g\w
X\not\simeq \w X$ for every $g\in G\backslash\{1\}$.
\end{prop}
\noindent{\textbf{Proof:}} We have $\ts{dim}_k\ts{End}_{\dba}(X,X)=1$ because $X$ is
indecomposable and $\dba(X,X[i])=0$ for every $i\neq 0$. On the other
hand, the spaces $\bigoplus\limits_{g\in G}\dbc(\,^g\w X,\w X)$ and
$\ts{End}_{\dba}(X,X)$ are isomorphic. So $\dbc(\,^g\w X,\w X)=0$ and,
therefore, $^g\w X\not\simeq X$ if $g\neq 1$.\sq

\paragraph{Proof of assertion ($H_3$)}$\ $

If $\psi\colon A\to A$ is an automorphism (and therefore a
Galois covering with trivial group), then
$\psi_{\lambda}\colon\ts{mod}\, A\to\ts{mod}\, A$ is an
equivalence. It thus induces a triangle equivalence
$\psi_{\lambda}\colon\dba\to\dba$.
\begin{prop}
\label{prop3.6}
Let $A$ be as in \ref{cor3.4},
$\psi\colon A\xrightarrow{\sim}A$ an automorphism such that
$\psi(x)=x$ for every $x\in \ts{ob}(A)$ and $T\in\dba$ a tilting
complex. Then $\psi_{\lambda}X\simeq X$ in $\dba$
for every indecomposable summand $X$ of $T$.
\end{prop}
\noindent{\textbf{Proof:}} Since $\psi(x)=x$ for every $x\in\ts{ob}(A)$,
we have the following fact:
\begin{enumerate}
\item[($i$)]
The conclusion of the proposition holds true if
$X$ is an indecomposable projective $A$-module.
\end{enumerate}
 Recall that
$\Theta\colon \dbh\to\dba$ is a triangle equivalence. We keep the
notations $\c T$ and $\c T'$ introduced in the proof of \ref{prop3.1}.
We prove the proposition for any $T\in\c T'$. By construction of
$\Theta$, we have:
\begin{enumerate}
\item[($ii$)] $\Theta$ induces a bijection $\Theta\colon\c T\to\c
  T'$. Under this bijection, tilting complexes in $\dbh$ correspond to
  tilting complexes in $\dba$.
\end{enumerate}
Since $\psi_{\lambda}\colon\dba\to\dba$ is an equivalence, we also have:
\begin{enumerate}
\item[($iii$)] Let $T,T'\in\c T$ be such that $T'$ is obtained from $T$
  by a transformation of the first kind. Then the proposition holds true
  for $T$ if and only if it does for $T'$.
\end{enumerate}
Now assume that $T,T'\in\c T'$ are such that $\Theta^{-1}(T')$ is
obtained from $\Theta^{-1}(T)$ by a transformation of the second
kind. We prove that the proposition holds true for $T$ if and
only if it does for $T'$. There exist
$X,Y\in\dba$ indecomposables and $\overline{T}\in\dba$ such
that $T=X\oplus\overline{T}$ and $T'=Y\oplus\overline{T}$. Also, there
exists a triangle in $\dba$ of one the two following forms:
\begin{enumerate}
\item $X\to M\to
  Y\to X[1]$  where $X\to M$ is a
  left minimal
  $\ts{add}(\overline{T})$-approximation and $M\to Y$ is a right minimal
   $\ts{add}(\overline{T})$-approximation.
\item $Y\to M\to
  X\to Y[1]$  where $Y\to M$ is a
  left minimal
  $\ts{add}(\overline{T})$-approximation and $M\to X$ is a right minimal
   $\ts{add}(\overline{T})$-approximation.
\end{enumerate}
Assume that the proposition holds for $T$ and that there is a triangle
$X\to M\to Y\to X[1]$ (the other cases are dealt with using similar
arguments). We only need to prove that $\psi_{\lambda}Y\simeq Y$. Apply
$\psi_{\lambda}$ to the triangle $X\to M\to Y\to X[1]$. Since
$\psi_{\lambda}$ is an equivalence and  the proposition holds
true for $T$, there exists a triangle $X\to M\to \psi_{\lambda}Y\to X[1]$ in $\dba$
where $X\to M$ is a left minimal
$\ts{add}(\overline{T})$-approximation. Therefore $\psi_{\lambda}Y\simeq
Y$ in $\dba$. So we proved that:
\begin{enumerate}
\item[($iv$)] If $T,T'\in\c T'$ are such that $\Theta^{-1}(T')$ is
  obtained from $\Theta^{-1}(T)$ by a transformation of
  the second kind, then the proposition holds true for $T$ if and only
  if it does for $T'$.
\end{enumerate}
As in the proof of \ref{prop3.1}, the conclusion follows
from ($i$), ($ii$), ($iii$), ($iv$) and \ref{prop2.4}.\sq

\paragraph{Application: proof of Theorem~\ref{thm0.3}}$\ $

As an application of the preceding results of
the section, we prove Theorem~\ref{thm0.3}. We  need the following
lemma. If $T=T_1\oplus\ldots\oplus T_n\in\ts{mod}\, A$ is an indecomposable
decomposition of a multiplicity-free module $T$, then $\ts{End}_A(T)$
is naturally a $k$-category,
equal to the full subcategory of $\ts{mod}\, A$ with objects $T_1,\ldots,T_n$.
\begin{lem}
  \label{lem3.7}
Let $A$ be a piecewise hereditary algebra of type $Q$.  Let $F\colon\c C\to A$ be a connected
  Galois covering with group $G$. Let $T\in\dba$ be a tilting complex, $B=\ts{End}_{\dba}(T)$
  and $T=T_1\bigoplus\ldots\bigoplus
  T_n$  an indecomposable decomposition. Let
  $\lambda_i\colon \l\w T_i\to T_i$ be an isomorphism where $\w
  T_i\in\dbc$ is indecomposable for every $i$. 
This defines the bounded complex of (not necessarily finite
dimensional) $\c C$-modules $\w T:=\bigoplus\limits_{i,g}\,^g\w T_i$,
where the sum runs over $g\in G$ and $i\in\{1,\ldots,n\}$. Let
$\c C'$ be the full subcategory of $\dbc$ with objects
the complexes $^g\w T_i$ (for $g\in G, i\in\{1,\ldots,n\}$).
Then the triangle
  functor $\l\colon\dbc\to \dba$ induces a connected Galois covering
  with group $G$:
\begin{equation}
\begin{array}{crclc}
  F_{\w T,\lambda}\colon&\c C'&\to&
B&\\
&^g\w T_i & \mapsto & T_i&\\
&^g\w T_i\xrightarrow{u} \,^h\w T_j & \mapsto  & T_i\xrightarrow{\
  \lambda_j\ \l(u)\ \lambda_i^{-1}} T_j& .
\end{array}\notag
\end{equation}
The complex $\w T$ is naturally a bounded complex of $\c C'-\c
C$-bimodules: As a functor from $\ts{End}_{\dbc}(\w T)\times\c
C^{op}$, it assigns the vector space $^g\w T_i(x)$ to the pair of
objects $(\,^g\w T_i,x)$. The total derived functor:
\begin{equation}
  -\stackrel[\c C']{\mathbb{L}}{\otimes}\w T\colon \c
  D^b(\ts{mod}\,\c C')\to \dbc\notag
\end{equation}
is  a $G$-equivariant triangle equivalence.
Finally, if $T$ is a tilting $A$-module and all the objects $\w T_i$
are $\c C$-modules (see \ref{cor3.4}), then:
\begin{enumerate}
\item[(a)] $\ts{Ext}^1_{\c C}(\,^g\widetilde{T}_i,\,^h\widetilde{T}_j)=0$ for every
  $i,j\in\{1,\ldots,n\}$ and $g,h\in G$.
\item[(b)] $\ts{pd}_{\c C}(\,^g\widetilde{T}_i)\leqslant 1$ for every $i,g$.
\item[(c)] If $P\in\ts{mod}\, \c C$ is an indecomposable projective $\c
  C$-module, then there exists an exact sequence $0\to P\to
  T^{(1)}\to T^{(2)}\to 0$ in $\ts{mod}\,\c C$ where $T^{(1)},T^{(2)}\in add(\{\,^g\widetilde{T}_i\ |\ i\in\{1,\ldots,n\},\ g\in G\})$.
  \end{enumerate}
\end{lem}
\noindent{\textbf{Proof:}} 
 By \ref{prop1.1}, the functor $F_{\w T,\lambda}$
 is a
well-defined Galois covering. By \ref{prop3.5}, we know that that
$\c C'$ is a locally bounded $k$-category (see
 \cite[2.1]{lemeur5}, for more details on the construction of
 $F_{\w T,\lambda}$). We prove that $\c C'$ is connected.
 By definition of $\w T$, we have $^g\w T=\w T$ for every $g\in
 G$. Hence the functor
$-\stackrel[\c C']{\mathbb{L}}{\otimes}\w T$ is $G$-equivariant. On
the other hand, $-\stackrel[\c C']{\mathbb{L}}{\otimes}\w T$ is a
triangle equivalence. Indeed, by
\ref{prop3.1}, (d), and by classical
arguments  on derived
equivalences (see \cite[III.2]{happel_book}, for instance), this
functor is full and faithful. Moreover its image
contains the complexes $^g\w T_i$ (for $g\in 
G$ and $i\in\{1,\ldots,n\}$). So \ref{prop3.1}, (b),
implies that this functor is dense and, therefore, a triangle equivalence $\c
D^b(\ts{mod}\, \c C')\to\dbc$. In
particular, $\c C'$ is connected.

Now we assume that $T$ is a tilting $A$-module. Assertion (a) follows from
\ref{prop3.1},
(d). Assertion (b) follows from the fact that
$\ts{pd}_A(T)\leqslant 1$ and $\l\colon\ts{mod}\, \c C\to \ts{mod}\, A$ is
exact. We prove assertion (c). Let $P\in\ts{mod}\, \c C$ be indecomposable
projective. Since $\l P$ is projective, there
exists an exact sequence $0\to \l P\to X\to Y\to 0$ in $\ts{mod}\, A$ with $X,Y\in
\ts{add}(T)$. By \ref{lem3.2},  the triangle $\l P\to X\to
Y\to \l P[1]$ is isomorphic to the image under $\l$ of a triangle $P\to
X'\to Y'\to P[1]$ where $X',Y'\in\ts{add}(\{\,^g\w T_i\ |\ g\in
G\,,i\in\{1,\ldots,n\}\})$. Since $\l$ is exact, the sequence $0\to
P\to X'\to Y'\to 0$ is an exact sequence in $\ts{mod}\,\c C$.
\sq

\begin{rem}
  \label{rem3.8}
Keep the hypotheses and notations of the preceding lemma.
  If $G$ is
finite and if $T$ is a tilting $A$-module, then the lemma expresses
that $\bigoplus\limits_{g,i}\,^g\w T_i$ is a tilting $\c C$-module.
\end{rem}

Now we can prove Theorem~\ref{thm0.3} which was
stated in the introduction.\\
\noindent{\textbf{Proof of Theorem~\ref{thm0.3}:}}
By \cite[Cor. 5.5]{happel_book}, there exists a sequence of
algebras:
\begin{equation}
A_0=kQ,A_1=\ts{End}_{A_0}(T^{(0)}),\ldots,A_l=\ts{End}_{A_{l-1}}(T^{(l-1)})=A\notag
\end{equation}
such that $T^{(i)}\in \ts{mod}\, A_{l-1}$ is tilting for every $i$. We
prove the theorem by induction on $l$. If
$l=0$, then $A=kQ$. For any connected Galois covering
$\c C\to A$ with group $G$ there exists a Galois covering of quivers
$Q'\to Q$ with group $G$ such that $\c C\simeq kQ'$ (see
\cite[Prop. 4.4]{lemeur2}). Assume that $l>0$ and the
conclusion of the theorem holds true for $A_{l-1}$. Let $\c C\to A$ be a
connected Galois covering with group $G$. Note that $T^{(l-1)}$ is a
tilting $A^{op}$-module. So the preceding lemma
yields a connected Galois covering $\c C'\to \ts{End}_{A^{op}}(T^{(l-1)})$
with group $G$ such that
$\c D^b(\ts{mod}\, \c C^{op})$ and $\c D^b(\ts{mod}\, \c C')$ are triangle
equivalent. 
 On the other hand,  $A_{l-1}\simeq
 \ts{End}_{A^{op}}(T^{(l-1)})^{op}$. Therefore the induction hypothesis implies
 that $\c D^b(\ts{mod}\,\c C'^{op})$ is triangle equivalent to $\c D^b(\ts{mod}\, kQ')$ where
 $Q'$ is a quiver such that there exists a Galois covering of quivers
 $Q'\to Q$ with group $G$.
 \sq

 \begin{rem}
   \label{rem3.9}
Let $A$ be a finite dimensional algebra endowed with a (non
necessarily free) $G$-action. Then:
\begin{enumerate}
\item[(a)]  If the $G$-action on $A$ is free, then the quotient algebra $A/G$ is well-defined.
  The proof of Theorem~\ref{thm0.3} shows that if $A/G$ is tilted
  (or, more generally,
  piecewise hereditary),
  then so is $A$. 
\item[(b)] It is proved in \cite[Thm. 3]{dionne_lanzilotta_smith} that if
  the order of $G$ is invertible in $k$ and if $A$ is piecewise
  hereditary, then so is the skew-group algebra $A[G]$. Recall that if
  $G$ acts freely on $A$, then the algebras $A[G]$ and $A/G$ are
  Morita equivalent (see \cite[Thm.2.8]{cibils_marcos}).
\end{enumerate}
 \end{rem}

\section{Correspondence between Galois coverings}
\label{sec:s5}$\ $

We still assume that $A$ is derived equivalent to a hereditary abelian
category $\c H$ such that $\overrightarrow{\c K}_{\c H}$ is connected.
Let
$T\in\dba$ be a tilting complex and $B=\ts{End}_{\dba}(T)$. In this
section, we construct a correspondence between the Galois coverings of
$A$ and those of $B$. This work has been done in \cite{lemeur5} in the
particular case where $T$ is a tilting $A$-module.
In order to compare the Galois coverings of $A$ and those of $B$,  it
is convenient to use the notion of
equivalent Galois covering. Given two Galois coverings $F\colon\c C\to
A$ and $F'\colon\c C'\to A$, we say that \textit{$F$ and $F'$ are
  equivalent} if there exists a commutative diagram:
\begin{equation}
  \xymatrix{
\c C \ar[r]^{\sim} \ar[d]_F & \c C'\ar[d]^{F'}\\
A \ar[r]_{\varphi}^{\sim} & A&
}\notag
\end{equation}
where the horizontal arrows are isomorphisms and $\varphi\colon
A\to A$ is an automorphism such that $\varphi(x)=x$ for every $x\in
\ts{ob}(A)$.

\paragraph{Equivalence classes of Galois coverings of $A$ associated
  to equivalence classes of Galois coverings of $B$}$\ $

In \ref{lem3.7}, we have associated a Galois covering $F_{\w
  T,\lambda}$ of $B$
to any Galois covering of $A$ and to any data consisting of
isomorphisms $(\lambda_i\colon \l\w T_i\to T_i)_{i=1,\ldots,n}$ in
$\dbc$. The following
lemma shows that different choices for these data give rise to
equivalent Galois coverings.
\begin{lem}{\cite[§ 2]{lemeur5}}
  \label{lem5.4}
Let $F\colon\c C\to A$ be a connected Galois covering with group
$G$. Let $T\in\dba$ be a tilting complex and $T=T_1\oplus\ldots\oplus
T_n$ be an indecomposable decomposition.
\begin{enumerate}
\item[(a)] Let $(\lambda_i\colon \l\w T_i\to T_i)_{i=1,\ldots,n}$ and
  $(\mu_i\colon \l\widehat{T}_i\to T_i)_{i=1,\ldots,n}$ be
  isomorphisms in $\dba$ defining the Galois coverings $F_{\w
    T,\lambda}\colon \ts{End}_{\dbc}(\w T)\to \ts{End}_{\dba}(T)$ and
  $F_{\widehat{T},\mu}\colon\ts{End}_{\dbc}(\widehat{T})\to
  \ts{End}_{\dba}(T)$, respectively.
Then  $F_{\w T,\lambda}$ and $F_{\widehat{T},\mu}$
  are equivalent. We write $[F]_T$ for the
  corresponding equivalence class of Galois coverings of $\ts{End}_{\dba}(T)$.
\item[(b)] Let $F'\colon\c C'\to A$ be a connected Galois covering with
  group $G$ and equivalent to $F$. Then the equivalence classes $[F]_T$
  and $[F']_T$ coincide.
\end{enumerate}
\end{lem}
\noindent{\textbf{Proof:}} In the case of tilting modules, (a)
and (b) were proved in \cite[Lem. 2.4]{lemeur5}
and \cite[Lem. 2.5]{lemeur5}, respectively. Using \ref{lem1.3}
and \ref{prop3.6}, it is
easily checked that the same arguments apply to prove (a)
and (b) in the present situation.\sq

In the sequel, we keep the notation $[F]_T$ introduced in \ref{lem5.4}. 

\paragraph{Galois coverings of $A$ induced by Galois coverings of
  $B$}$\ $

We now  express any Galois covering of $A$ as
induced by a Galois covering of $B$ as in \ref{lem3.7}. The
tilting complex $T$ is naturally a complex of $B-A$-bimodules. Also,
it defines a triangle equivalence:
\begin{equation}
  -\stackrel[B]{\mathbb{L}}{\otimes}X\colon \dbb\to\dba\ .\notag
\end{equation}
Fix a connected Galois covering $F\colon\c C\to A$ with
group $G$, an indecomposable decomposition
$T=T_1\oplus\ldots\oplus T_n$ and isomorphisms $(\mu_i\colon \l\w
T_i\xrightarrow{\sim} T_i)_{i=1,\ldots,n}$. According to
\ref{lem3.7}, these data define the Galois covering $F_{\w
  T,\mu}\colon \ts{End}_{\dbc}(\w T)\to B$ which we denote
by $F'\colon \c C'\to B$ for simplicity.
\begin{lem}
  \label{lem5.5}
The following diagram commutes up to an isomorphism of functors.
\begin{equation}
  \xymatrix{
\c D^b(\ts{mod}\, \c C') \ar[d]_{F'_{\lambda}} \ar[rr]^{-\stackrel[\c
  C']{\mathbb{L}}{\otimes}\w T} && \dbc \ar[d]^{\l}\\
\dbb \ar[rr]_{-\stackrel[B]{\mathbb{L}}{\otimes}T} && \dba &,
}\notag
\end{equation}
\end{lem}
\noindent{\textbf{Proof:}} Recall that $\l\colon \ts{mod}\,\c C\to\ts{mod}\, A$
(or $F'_{\lambda}\colon\ts{mod}\,\c C'\to\ts{mod}\, B$) is
equal to $-\stackrel[\c C]{}{\otimes}A$ (or to $-\stackrel[\c
C']{}{\otimes}B$, respectively). Since these two functors are exact
and map projective modules to projective modules and  the
horizontal arrows of the diagram are triangle equivalences (see
\ref{lem5.4}), we deduce that:
\begin{enumerate}
\item The composition $\c D^b(\ts{mod}\, \c C')\xrightarrow{F'_{\lambda}}\c
  D^b(B)\xrightarrow{-\stackrel[B]{\mathbb{L}}{\otimes}T}\dba$ is
  isomorphic to $-\stackrel[\c
  C']{\mathbb{L}}{\otimes}\left(B\stackrel[B]{}{\otimes}T\right)$.
\item The composition $\c D^b(\ts{mod}\, \c
  C')\xrightarrow{-\stackrel[\c C']{\mathbb{L}}{\otimes}\w T}\dbc\xrightarrow{\l}\dba$ is isomorphic to
  $-\stackrel[\c C']{\mathbb{L}}{\otimes}\left(\w T\stackrel[\c C]{}{\otimes}A\right)$.
\end{enumerate}
On the other hand, the isomorphisms
$\mu_i\colon\l\w T_i\xrightarrow{\sim} T_i$ (for $i\in\{1,\ldots,n\}$)
define an isomorphism $B\stackrel[B]{}{\otimes}T\xrightarrow{\sim}\w T\stackrel[\c
C]{}{\otimes}A$ of $\c C'-A$-bimodules. This proves that the diagram
commutes up to an isomorphism of functors.\sq

Since $-\stackrel[B]{\mathbb{L}}{\otimes}T$ is an equivalence, there
exists an isomorphism $\varphi_x\colon
X_x\stackrel[B]{\mathbb{L}}{\otimes}T\to A(-,x)$ in $\dba$ with
$X_x\in\dbb$ for every $x\in\ts{ob}(A)$. In particular,
$\bigoplus\limits_{x\in\ts{ob}(A)}X_x$ is an indecomposable decomposition
of a tilting complex in $\dbb$. Then by the preceding section,
there exists an isomorphism $\nu_x\colon F'_{\lambda}\w
X_x\xrightarrow{\sim} X_x$ in $\dbb$ with $\w X_x\in\c D^b(\ts{mod}\, \c
C')$ for every $x\in\ts{ob}(A)$. By \ref{lem3.7}, the datum
$(\nu_x)_{x\in\ts{ob}(A)}$ defines a connected Galois covering with group $G$:
\begin{equation}
  \begin{array}{crclc}
    F'_{\w X,\nu}\colon&\ts{End}_{\c D^b(\ts{mod}\, \c C')}(\w X) & \to &\ts{End}_{\dbb}(X)&\\
&^g\w X_x & \mapsto & X_x&\\
&^g\w X_x\xrightarrow{u}\,^h\w X_y & \mapsto &X_x\xrightarrow{\nu_{y}\
  F'_{\lambda}(u)\ \nu_x^{-1}}X_y&.
  \end{array}
  \notag
\end{equation}
On the other hand, the isomorphisms $\varphi_x$ (for $x\in \ts{ob}(A)$) define
the following isomorphism of $k$-categories:
\begin{equation}
  \begin{array}{crclc}
    \rho_{X,\varphi}\colon& \ts{End}_{\dbb}(X) & \to & A&\\
    &X_x & \mapsto & x&\\
&X_x\xrightarrow{u}X_y & \mapsto & \left(\varphi_y
  \circ(u\stackrel[]{}{\otimes}T)\circ \varphi_x^{-1}\right)(\ts{Id}_x)\in
A(x,y)&.
  \end{array}\notag
\end{equation}
Thus, we have a connected Galois
covering $\rho_{X,\varphi}\circ F'_{\w X,\nu}\colon\ts{End}_{\c D^b(\ts{mod}\,
  \c C')}(\w X)\to A$ with group $G$ which we denote by $F''$.
The following lemma relates $F$ and $F''$.
\begin{lem}
  \label{lem5.6}
The Galois coverings $F$ and $F''$ are equivalent.
\end{lem}
\noindent{\textbf{Proof:}} We need to construct a commutative diagram:
\begin{equation}
  \xymatrix{
\ts{End}_{\c D^b(\ts{mod}\, \c C')}(\w X) \ar[rr]^{\sim} \ar[d]_{F''} && \c
C\ar[d]^F\\
A \ar[rr]^{\sim} && A
}\tag{$D$}
\end{equation}
where the horizontal arrows are isomorphisms and  the bottom
horizontal isomorphism extends the identity map on objects.
For this purpose, we proceed in two steps. 

\textbf{Step $1$:} We express $F$ as a
functor between subcategories of $\dbc$ and $\dba$. Given $x\in \ts{ob}(\c
C)$, the $A$-module $\l(\c C(-,x))$ does depend
only on $F(x)$ (and not on $x$) because $\l$ is $G$-invariant. Besides,
there is a canonical isomorphism $\iota_x\colon \l(\c
C(-,x))\xrightarrow{\sim}A(-,F(x))$ of $A$-modules induced by $F$: If $y\in
\ts{ob}(A)$, then $\left(\l(\c
  C(-,x))\right)(y)=\bigoplus\limits_{Fy=F(y)}\c C(y',x)$ and an
element 
$(u_{y'})_{y'}$ of this vector space is mapped by $\iota_x$ to
$\sum\limits_{y'}F(u_{y'})\in A(F(y),F(x))$.
Clearly, this isomorphism does depend only on $F(x)$
(and not on $x$) whence the notation $\iota_x$. Now, let $\c P_A$ and $\c P_{\c C}$ be the full subcategories of
$\dba$ and $\dbc$ with object sets $\{A(-,x)\ |\ x\in \ts{ob}(A)\}$
and $\{\c C(-,x)\ |\ x\in \ts{ob}(\c C)\}$, respectively. Hence we have a
commutative diagram:
\begin{equation}
  \xymatrix{
\c C \ar[r]^{\sim} \ar[d]_F & \c P_{\c C} \ar[d] \\
A \ar[r]^{\sim} & \c P_A&
}\tag{$D_1$}
\end{equation}
where the unlabelled functors are as follows:
\begin{enumerate}
\item The functor $\c C\to \c P_{\c C}$ is the following isomorphism:
  \begin{equation}
    \begin{array}{rclc}
      \c C &\to & \c P_{\c C}&\\
      x\in \ts{ob}(\c C) & \mapsto & \c C(-,x)&\\
u\in\,\c C(x,y) & \mapsto & \c C(-,u)\colon \c C(-,x)\to\c C(-,y)&.
    \end{array}\notag
  \end{equation}
\item The functor $A\to \c P_{A}$ is the following isomorphism:
  \begin{equation}
    \begin{array}{rclc}
      A &\to & \c P_{A}&\\
      x\in \ts{ob}(A) & \mapsto & A(-,x)&\\
u\in\,A(x,y) & \mapsto & A(-,u)\colon A(-,x)\to A(-,y)&.
    \end{array}\notag
  \end{equation}
\item The functor $\c P_{\c C}\to \c P_A$ is as follows:
  \begin{equation}
    \begin{array}{rclc}
      \c P_{\c C} & \to &\c P_A& \\
      \c C(-,x) & \mapsto &A(-,F(x))&\\
      \c C(-,x)\xrightarrow{u}\c C(-,y) & \mapsto &
      A(-,F(x))\xrightarrow{\iota_y\circ \l(u)\circ \iota_x^{-1}}A(-,F(y))&.
    \end{array}\tag{$i$}
  \end{equation}
\end{enumerate}
In particular, $\c P_{\c C}\to \c P_A$ is a Galois covering with group $G$.

\textbf{Step 2:} We now relate $F''$ to the
Galois covering $\c P_{\c C}\to \c P_A$.
We first construct an isomorphism $\ts{End}_{\c D^b(\ts{mod}\, \c C')}(\w
X)\xrightarrow{\sim}\c P_{\c C}$. Let $\Theta\colon
\l'(-)\stackrel[B]{\mathbb{L}}{\otimes}T\xrightarrow{\sim}
  \l(-\stackrel[\c C']{\mathbb{L}}{\otimes}\w T)$ be an isomorphism of functors
  (see \ref{lem5.5}). Let $x\in
\ts{ob}(A)$. So we have a composition of isomorphisms in $\dba$:
\begin{equation}
  \l\left(\w X_x\stackrel[\c C']{\mathbb{L}}{\otimes}\w
    T\right)\xrightarrow{\Theta_{\w X_x}^{-1}}F'_{\lambda}\w
    X_x\stackrel[B]{\mathbb{L}}{\otimes}T \xrightarrow{\nu_x\otimes
      T}X_x\stackrel[B]{\mathbb{L}}{\otimes}T\xrightarrow{\varphi_x}A(-,x)\ .\notag
\end{equation}
Therefore, by \ref{lem1.3}, there exists  an isomorphism
$\psi_x\colon \w X_x\stackrel[\c C']{\mathbb{L}}{\otimes}\w
    T\xrightarrow{\sim}\c C(-,L(x))$ in $\dbc$ with $L(x)\in
    F^{-1}(x)$. We deduce that the following is an
    isomorphism of $k$-categories because $-\stackrel[\c
    C']{\mathbb{L}}{\otimes}\w T$ is a $G$-equivariant functor (see \ref{lem3.7}):
    \begin{equation}
      \begin{array}{rclc}
         \ts{End}_{\c D^b(\ts{mod}\, \c C')}(\w X) & \to & \c P_{\c C} &\\
       ^g\w X_x & \mapsto & \c C(-,gL(x))&\\
       ^g\w X_x\xrightarrow{u}\,^h\w X_y & \mapsto &
        \c C(-,gL(x))\xrightarrow{^h\psi_y\circ (u\otimes \w T)\circ
          (\,^g\psi_x)^{-1}}\c C(-,hL(y))& .
      \end{array}\tag{$ii$}
    \end{equation}
We now construct another isomorphism between $A$ and $\c P_A$.
We have the following composition of isomorphisms in $\dba$
which we denote by $\alpha_x$:
\begin{equation}
\alpha_x\colon   A(-,x)
  \xrightarrow{\varphi_x^{-1}}X_x\stackrel[B]{\mathbb{L}}{\otimes}T
  \xrightarrow{(\nu_x\otimes T)^{-1}} F'_{\lambda}\w
  X_x\stackrel[B]{\mathbb{L}}{\otimes} T \xrightarrow{\Theta_{\w X_x}}
  \l\left(\w 
    X_x\stackrel[\c C']{\mathbb{L}}{\otimes}\w T\right)
  \xrightarrow{\l(\psi_x)} \l(\c C(-,L(x))\xrightarrow{\iota_x}A(-,x).\notag
\end{equation}
Note that $\alpha_x\colon A(-,x)\xrightarrow{\sim}A(-,x)$ is
necessarily equal to the multiplication by a scalar in $k^*$ because $A(-,x)$ is
an indecomposable projective $A$-module and $A$ is piecewise
hereditary. Therefore we have an isomorphism of categories:
\begin{equation}
  \begin{array}{rclc}
    A & \to & \c P_A & \\
    x & \mapsto & A(-,x)&\\
    u\in A(x,y) & \mapsto & \alpha_y\circ A(-,u)\circ \alpha_x^{-1}&.
  \end{array}\tag{$iii$}
\end{equation}
Hence the horizontal arrows of the following diagram are
isomorphisms:
\begin{equation}
  \xymatrix{
\ts{End}_{\c D^b(\ts{mod}\, \c C')}(\w X) \ar[d]_{F''} \ar[rr]^{\text{\normalfont ($ii$)}}&& \c
P_{\c C} \ar[d]^{\text{\normalfont ($i$)}} \\
A  \ar[rr]_{\text{\normalfont ($iii$)}} && \c P_A & .
}\tag{$D_2$}
\end{equation}
We claim that this diagram commutes. The commutativity is clearly satisfied
on objects. Let $u\colon \,^g\w X_x\to\,^h\w X_y$ be a morphism in
$\ts{End}_{\c D^b(\ts{mod}\, \c C')}(\w X)$. Denote by $u_1\colon A(-,x)\to
A(-,y)$ the image of $u$ under the
composition of ($i$) and ($ii$). Then:
\begin{equation}
\begin{array}{rcllc}
  u_1&=&\iota_y\circ\l\left(^h\psi_y\circ (u\otimes \w T)\circ
          (\,^g\psi_x)^{-1}\right)\circ \iota_x^{-1}&&\\
&=&\iota_y\circ\l(\psi_y)\circ\l\left(u\otimes \w
  T\right)\circ\left(\l(\psi_x)\right)^{-1}\circ\iota_x^{-1}&\
\text{because $\l$ is $G$-invariant,}&\\
&=&\iota_y\circ\l(\psi_y)\circ\Theta_{\w X_y}\circ \left(F'_{\lambda}(u)\otimes
  T\right)\circ \Theta_{\w
  X_x}^{-1}\circ\left(\l(\psi_x)\right)^{-1}\circ\iota_x^{-1}&\
\text{by definition of $\Theta$,}&\\
&=&\alpha_y\circ \varphi_y\circ \left(\nu_y\otimes T\right)\circ \left(F'_{\lambda}(u)\otimes
  T\right) \circ \left(\nu_x\otimes
  T\right)^{-1}\circ \varphi_x^{-1}\circ \alpha_x^{-1}&\ \text{by
  definition of $\alpha_x$ and $\alpha_y$,}&\\
&=&\alpha_y\circ \varphi_y\circ \left((\nu_y\circ F'_{\lambda}(u)\circ
  \nu_x^{-1})\otimes T\right)\circ \varphi_x^{-1}\circ
\alpha_x^{-1}&&.\\
\end{array}\notag
\end{equation}
Now, let $u_2\in A(x,y)$ be the image of $u$ under $F''$, that is
$u_2=\left(\varphi_y\circ \left((\nu_y\circ F'_{\lambda}(u)\circ
  \nu_x^{-1})\otimes T\right)\circ
\varphi_x^{-1}\right)(\ts{Id}_x)$. Therefore $A(-,u_2)$ is equal to the
morphism $\varphi_y\circ \left((\nu_y\circ F'_{\lambda}(u)\circ
  \nu_x^{-1})\otimes T\right)\circ
\varphi_x^{-1}\colon A(-,x)\to A(-,y)$ of $\c P_A$. In particular,
the image of $u_2$ under ($iii$) coincides with $u_1$. Therefore
($D_2$) is commutative. Since ($D_1$) also commutes, we deduce
that so does ($D$). Thus, $F$ and $F''$ are equivalent.\sq

\paragraph{Correspondence between the Galois coverings of $A$ and
  those of $B$}$\ $

\begin{prop}
  \label{prop5.7}
Let $A$ be an algebra derived equivalent to a hereditary abelian
category $\c H$ such that $\overrightarrow{\c K}_{\c H}$ is connected.
Let $T\in\dba$ be a tilting complex, $B=\ts{End}_{\dba}(T)$ and $G$
 a group. With the
notations of \ref{lem5.4}, the map $[F]\mapsto [F]_T$ is a
well-defined bijection from the set of equivalence classes of
connected Galois coverings with group $G$ of $A$ to the set of
equivalence classes of Galois coverings with group $G$ of $B$.
\end{prop}
\noindent{\textbf{Proof:}}  Let $\ts{Gal}_A(G)$ be the set of equivalence
classes of connected Galois coverings with group $G$ of $A$. By
\ref{lem5.4}, there is a well-defined map:
\begin{equation}
  \begin{array}{crclc}
    \gamma_A\colon & \ts{Gal}_A(G) & \to& \ts{Gal}_B(G)\\
    &[F] & \mapsto & [F]_T &.
  \end{array}\notag
\end{equation}
We keep the notations $X_x,\varphi_x$ (for $x\in \ts{ob}(A)$)  introduced
 after the proof of \ref{lem5.5}. Then we also have a
well-defined map:
\begin{equation}
  \begin{array}{crclc}
    \gamma_B\colon & \ts{Gal}_B(G) & \to &\ts{Gal}_{\ts{End}_{\dbb}(X)}(G)\\
    &[F] & \mapsto & [F]_X 
  \end{array}\notag
\end{equation}
 By \ref{lem5.6}, we know that $\gamma_A$ is injective and
$\gamma_B$ is surjective. Therefore $\gamma_A$ is bijective because
$A,T$ and $B,X$ play symmetric rôles.\sq

By \ref{prop5.7}, we have some information on the existence of a
universal cover. Indeed, we have the following result.
\begin{prop}
  \label{prop5.3}
Let $A$ be as in \ref{prop5.7} and
  $T\in\dba$ a tilting complex. Assume that $A$ admits a universal cover
$\widetilde{F}\colon\widetilde{\c C}\to A$. Then
$\ts{End}_{\dba}(T)$ admits a universal cover with group isomorphic to the one
of $\widetilde{F}$.
\end{prop}
\noindent{\textbf{Proof:}} Fix an indecomposable decomposition $T=T_1\oplus\ldots \oplus
T_n$. Let $B=\ts{End}_{\dba}(T)$. So $B$ is the full subcategory of $\dba$ with objects
$T_1,\ldots,T_n$. Let $x_0\in \ts{ob}(A)$ be a base-point for the category $\ts{Gal}(A,x_0)$
of pointed Galois coverings of $A$. We  construct a (full and faithful) functor $\widetilde{F}^{\to}\to
\ts{Gal}(B,T_1)$. Recall that
$\widetilde{F}^{\to}$ was defined in Section~\ref{sec:s1} and there
is at most one morphism between two pointed Galois coverings. We need
the following
data:
\begin{enumerate}
\item For every $i\in\{1,\ldots,n\}$, let $\widetilde{T}_i\in
  \ts{mod}\,\widetilde{\c C}$ be such that
  $\widetilde{F}_{\lambda}\widetilde{T_i}\simeq T_i$. Therefore the
  $k$-categories $B=\ts{End}_{\dba}(T)$ and
  $\ts{End}_{\dba}(\bigoplus\limits_{i=1}^n\l\w T_i)$ are
  isomorphic. For simplicity, we assume that
  $\widetilde{F}_{\lambda}\widetilde{T_i}=T_i$ for every $i$.
\item If $F\in \widetilde{F}^{\to}$, there exists a unique
  morphism $p\colon \widetilde{F}\to F$ in $\ts{Gal}(A,x_0)$. Since $p$ is
  a Galois covering (see \cite[Prop. 3.4]{lemeur2}), we set
  $T^F_i=p_{\lambda}\widetilde{T}_i$ for every $i$.
\end{enumerate}
Then:
\begin{enumerate}
\item[($i$)] We have
  $T_i=F_{\lambda}(T^F_i)$ for every $i\in\{1,\ldots,n\}$ and
  $F\in\w F^{\to}$. Indeed, there exists a unique morphism $p\colon \w
  F\to F$, so that $\w F_{\lambda}=F_{\lambda}\circ p_{\lambda}$.
\item[($ii$)] Let $u\colon F\to F'$ be a morphism in
  $\widetilde{F}^{\to}$. Let $G$ be the group of $F$ and $G'$ 
  the group of $F'$. Then $u$ is a Galois covering
  (see \cite[Prop. 3.4]{lemeur2}). Let $p\colon \widetilde{F}\to F$ and
  $p'\colon \widetilde{F}\to F'$ be the unique morphisms in
  $\ts{Gal}(A,x_0)$ from $\w F$ to $F$ and from $\w F$ to $F'$,
  respectively. Then $p'=u\circ p$ because of \cite[Lem. 3.1]{lemeur2}
  and because we are dealing with pointed Galois coverings. Therefore
  $u_{\lambda}(T^F_i)=\,^{\sigma_u(g)}T_i^{F'}$ for every
  $i\in\{1,\ldots,n\}$ and every $g\in G$.
Here, $\sigma_u\colon G\to G'$ is the unique (surjective) morphism of
groups such that $u\circ g=\sigma_u(g)\circ u$ for every $g\in G$ (see \cite[Prop. 3.4]{lemeur2}).
\end{enumerate}

Now we can construct a functor $\widetilde{F}^{\to}\to
\ts{Gal}(B,T_1)$. Let $F\colon(\c C,x)\to (A,x_0)$ be in
$\widetilde{F}^{\to}$. Let $G$ be the group of $F$. By ($i$) and
\ref{lem3.7}, we
have a pointed Galois covering with group $G$ induced by
$\l\colon\dbc\to\dba$:
\begin{equation}
  \begin{array}{crclc}
    F_T\colon & \left(\ts{End}_{\dbc}(\bigoplus\limits_{g,i}\,^gT^F_i),\, T_1^F\right) & \to &
    (B,T_1)\\
& ^gT_i^F & \mapsto & T_i\\
& ^gT_i^F\xrightarrow{f}\,^hT_j^F & \mapsto & T_i\xrightarrow{\l(f)}T_j&.
  \end{array}\notag
\end{equation}
So
$[F_T]=[F]_T$.
 Thus, we have associated a pointed Galois
covering with group $G$ of $B$ to any pointed Galois covering with
group $G$ of $A$. We now associate a morphism of pointed Galois
coverings of $B$ to any morphism of pointed Galois coverings of $A$. Let $u\colon F\to
F'$ be a morphism in $\widetilde{F}^{\to}$ where $F\colon(\c C,x)\to (A,x_0)$ and $F'\colon (\c
C',x')\to (A,x_0)$ have groups $G$ and
$G'$, respectively. By ($ii$), we
have a well-defined $k$-linear functor induced by $u_{\lambda}\colon
\dbc\to \c D^b(\ts{mod}\, \c C')$:
\begin{equation}
  \begin{array}{crclc}
    u_T\colon &
    \left(\ts{End}_{\dbc}\left(\bigoplus\limits_{g,i}\,^gT_i^F\right),\,
      T_1^F\right) & \to &\left(\ts{End}_{\c
        D^b(\ts{mod}\,\c C')}\left(\bigoplus\limits_{g',i}\,^{g'}T_i^{F'}\right),\,
      T_1^{F'}\right)\\
    & ^gT_i^F & \mapsto  & u_{\lambda}(\,^gT_i^F)=\,^{\sigma_u(g)}T_i^{F'}\\
    & ^gT_i^F\xrightarrow{f}\,^hT_j^F & \mapsto
    &\,^{\sigma_u(g)}T_i^{F'}\xrightarrow{u_{\lambda}(f)}\,^{\sigma_u(h)}T_j^{F'}&.
  \end{array}\notag
\end{equation}
The equality $u_{\lambda}(\,^gT_i^F)=\,^{\sigma_u(g)}T_i^{F'}$
follows from the equality $u\circ g=\sigma_u(g)\circ u$. Also,
$u_{\lambda}(T_1^F)=T_1^{F'}$.  Since  $F'\circ u=F$ and $F_T,F'_T$ and  
$u_T$ are defined as restrictions of $F_{\lambda}$,
$F'_{\lambda}$  and
$u_{\lambda}$ respectively, 
$u_T\colon F_T\to F'_T$ is a morphism in $\ts{Gal}(B,T_1)$. Thus, to any
morphism in $\widetilde{F}^{\to}$, we have associated a morphism in
$\ts{Gal}(B,T_1)$. We let the
reader check that the following is a functor:
\begin{equation}
  \begin{array}{crcl}
    \Psi \colon & \widetilde{F}^{\to} & \to & \ts{Gal}(B,T_1)\\
    & F & \mapsto & F_T\\
    & F\xrightarrow{u} F' & \mapsto & F_T\xrightarrow{u_T} F'_T
  \end{array}\notag
\end{equation}
Also, it is not difficult to prove that $\Psi$ is full and faithful,
although we shall not use this fact in the sequel. Remark that the
Galois covering $F_T$ lies in $\Psi(\widetilde{F})^{\to}$ for every
$F\in\widetilde{F}^{\to}$.\\

Now we can prove that $\Psi(\widetilde{F})$ is a
universal cover for $B$. Let $F$ be a connected Galois covering of
$B$. By \ref{prop5.7}, there exists a connected Galois
covering $F'$ of $A$ such that $[F]=[F']_T$. Since $\widetilde{F}$ is
a universal cover of $A$, the Galois covering $F'$ of $A$ is
equivalent to some $F''\in \widetilde{F}^{\to}$, that is
$[F']=[F'']$. As noticed above, we have  $[F''_T]=[F'']_T$. Therefore
$[F]=[F']_T=[F'']_T=[F''_T]$, that
is, $F$ is equivalent to a Galois covering of $B$ lying in
$\Psi(\widetilde{F})^{\to}$. So $\Psi(\widetilde{F})$ is a
universal Galois covering of $B$.\sq

\section{The main theorem and its corollary}
\label{sec:s6}$\ $

In this section, we prove Theorem~\ref{thm0.2}
and Corollary~\ref{cor0.1}. We assume that $A$ is a 
connected algebra derived equivalent to a hereditary abelian category
$\c H$ such that $\overrightarrow{\c K}_{\c H}$  is connected.
\paragraph{Two particular cases: paths algebras and squid algebras}$\
$

We first check that our main results hold for paths algebras and for
squid algebras.
\begin{lem}
  \label{lem6.1}
Assume that $A=kQ$ where $Q$ is a finite connected quiver with no
oriented cycle. Then Theorem~\ref{thm0.2} and Corollary~\ref{cor0.1}
hold true for $A$.
\end{lem}
\noindent{\textbf{Proof:}} Let $\w Q\to Q$ be the universal Galois
covering of quivers (see \cite{martinezvilla_delapena}). It follows
from \cite[Prop. 4.4]{lemeur2} that the induced Galois covering $k\w
Q\to kQ$ with group $\pi_1(Q)$ is a universal cover of $A$. Whence
Theorem~\ref{thm0.3}. On the other hand, $\ts{HH}^1(kQ)=0$ if
and only if $Q$ is a tree (see \cite{cibils_h1}). Whence
Corollary~\ref{cor0.1}.\sq

We now turn to the case of squid algebras. We refer the reader to \cite{ringel} for more details on squid
algebras. A squid algebra  over an algebraically closed field $k$
is defined by the following data: An integer $t\geqslant 2$, a
sequence $p=(p_1,\ldots,p_t)$ of non negative integers and a sequence
$\tau=(\tau_3,\ldots,\tau_t)$ of pairwise distinct non-zero elements of
$k$. With this data, the squid algebra $S(t,p,\tau)$ is the
$k$-algebra $kQ/I$ where $Q$ is the following quiver:
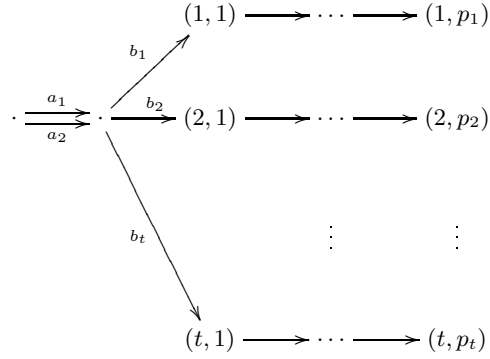
\begin{equation}
  \xymatrix{
&&(1,1)\ar@{->}[r]&\ldots \ar@{->}[r]& (1,p_1)\\
\cdot\ar@<2pt>@{->}[r]^{a_1}\ar@<-2pt>@{->}[r]_{a_2} & \cdot\ar@{->}[ru]^{b_1}
\ar@{->}[r]^{b_2}
\ar@{->}[rdd]_{b_t}&(2,1)\ar@{->}[r]&\ldots \ar@{->}[r]& (2,p_2)\\
&&&\vdots&\vdots\\
&&(t,1)\ar@{->}[r]&\ldots \ar@{->}[r]& (t,p_t)&
}\notag
\end{equation}
and $I$ is the ideal generated by the following relations:
\begin{equation}
  b_1a_1=b_2a_2=0,\ \ b_ia_2=\tau_i\, b_ia_1\ \ \text{for
    $i=3,\ldots,t$}\ .\notag
\end{equation}

Using Happel's long exact sequence (\cite{happel}), one can compute
$\ts{HH}^1(S(t,p,\tau))$:
\begin{equation}
  \operatorname{\ts{dim}_k\ \ts{HH}^1}(S(t,p,\tau))=\begin{cases}
1 & \text{if $t=2$}\\
0 & \text{if $t\geqslant 3$.}
\end{cases}\notag
\end{equation}
Following \cite{martinezvilla_delapena}, the bound quiver $(Q,I)$
defines a Galois covering $k\w Q/\w I\to kQ/I$ with group isomorphic
to $\mathbb{Z}$ if $t=2$ and with trivial group otherwise. One can
easily check that this Galois covering is universal in the sense of
Theorem~\ref{thm0.3}. The above considerations give the following.
\begin{lem}
  \label{lem6.2}
Let $A$ be a squid algebra. Then Theorem~\ref{thm0.2} and
Corollary~\ref{cor0.1} hold true for $A$.
\end{lem}

\paragraph{The general case}$\ $

Using \ref{prop5.3}, \ref{lem6.1} and \ref{lem6.2} we can prove the
 two main results of this text.\\
\noindent{\textbf{Proof of Theorem~\ref{thm0.2} and Corollary~\ref{cor0.1}:}} 
Assume first that $A$ is piecewise hereditary of type $Q$ where $Q$
is a finite connected quiver with no oriented cycle. So there exists
 a tilting complex $T\in\c D^b(\ts{mod}\, kQ)$ such that $A\simeq \ts{End}_{\c
  D^b(\ts{mod}\, kQ)}(T)$. By \ref{prop5.3} and \ref{lem6.1},
the algebra $A$ admits a universal Galois covering with group
isomorphic to the fundamental group of $Q$. In particular, $A$ is
simply connected if and only if $Q$ is a tree. On the other hand, $Q$
is tree if and only if $\ts{HH}^1(kQ)=0$ (by \cite{cibils_h1}) and
$\ts{HH}^1(kQ)\simeq \ts{HH}^1(A)$ (by \cite{keller}). Therefore $A$ is
simply connected if and only if $\ts{HH}^1(A)=0$, or, if and only if $Q$
is a tree.

Assume now that $A$ is not derived equivalent to $\c D^b(\ts{mod}\, kQ)$ for any
finite quiver $Q$. Then  \cite[Prop. 2.1,
Thm. 2.6]{happel_reiten} implies that there exists a squid algebra
$S=S(t,p,\tau)$ and a tilting complex $T\in\c D^b(\ts{mod}\, S)$ such that $A\simeq
\ts{End}_{\c D^b(\ts{mod}\, S)}(T)$. By \ref{prop5.3} and
\ref{lem6.2}, the algebra $A$ has a universal cover with
group isomorphic to the trivial group or to $\mathbb{Z}$ according to whether
$t\geqslant 3$ or $t=2$. In particular, $A$ is simply connected if and
only if $t=2$, that is, if and only if $\ts{HH}^1(S)=0$ (see
\ref{lem6.2}). Since $\ts{HH}^1(S)\simeq \ts{HH}^1(A)$
(by \cite{keller}), we deduce that $A$ is simply connected if and only if $\ts{HH}^1(A)=0$.
\sq

\bibliographystyle{plain}
\bibliography{biblio}

\begin{thebibliography}{10}

\bibitem{assem_marcos_delapena}
I.~Assem, E.~N. Marcos, and J.~A. de~La~Pe\~na.
\newblock The simple connectedness of a tame tilted algebra.
\newblock {\em J. Algebra}, 237(2):647--656, 2001.

\bibitem{assem_simson_skowronski}
I.~Assem, D.~Simson, and A.~Skowro\'nski.
\newblock {\em {Elements of the representation theory of associative algebras.
  Vol. 1: Techniques of representation theory.}}
\newblock {London Mathematical Society Student Texts 65. Cambridge: Cambridge
  University Press. ix, 458~p. }, 2006.

\bibitem{assem_skowronski}
I.~Assem and A.~Skowro\'nski.
\newblock On some classes of simply connected algebras.
\newblock {\em Proc. {L}ondon {M}ath. {S}oc.}, 56(3):417--450, 1988.

\bibitem{bongartz_gabriel}
K.~Bongartz and P.~Gabriel.
\newblock Covering spaces in representation theory.
\newblock {\em Invent. {M}ath.}, 65:331--378, 1982.

\bibitem{bmrrt}
A.~B. Buan, R.~Marsh, M.~Reineke, I.~Reiten, and G.~Todorov.
\newblock {Tilting theory and cluster combinatorics}.
\newblock {\em {Adv. Math.}}, 204(2):572--618, 2006.

\bibitem{buchweitz_liu}
R.-O. Buchweitz and S.~Liu.
\newblock {Hochschild cohomology and representation-finite algebras.}
\newblock {\em {Proc. London Math. Soc., Ser. III}}, 88(2):355--380, 2004.

\bibitem{ccs}
P.~Caldero, F.~Chapoton, and R.~Schiffler.
\newblock {Quivers with relations arising from clusters ($A_n$ case)}.
\newblock {\em Trans. Amer. Math. Soc.}, 2006.

\bibitem{cibils_h1}
C.~Cibils.
\newblock {On the Hochschild cohomology of finite dimensional algebras.}
\newblock {\em Comm. Algebra}, 16(3):645--649, 1988.

\bibitem{cibils_marcos}
C.~Cibils and E.~N. Marcos.
\newblock Skew categories, {G}alois coverings and smash-product of a
  $k$-category.
\newblock {\em Proc. {A}mer. {M}ath. {S}oc.}, 134(1):39--50, 2006.

\bibitem{dionne_lanzilotta_smith}
J.~Dionne, M.~Lanzilotta, and D.~Smith.
\newblock {Skew group algebras of piecewise hereditary algebras are piecewise
  hereditary}.
\newblock {\em {J. Pure Appl. Algebra}}, 213(2):241--249, 2009.

\bibitem{dowbor_skowronski}
P.~Dowbor and A.~Skowro\'nski.
\newblock Galois coverings of representation infinite algebras.
\newblock {\em {C}omment. {M}ath. {H}elv.}, 62:311--337, 1987.

\bibitem{fz}
S.~Fomin and A.~Zelevinsky.
\newblock {Cluster algebras I: Foundations}.
\newblock {\em Trans. Amer. Math. Soc.}, 353:497--529, 2002.

\bibitem{gabriel}
P.~Gabriel.
\newblock {The universal cover of a representation-finite algebra.}
\newblock {Representations of algebras, Proc. 3rd int. Conf., Puebla/Mex. 1980,
  Lecture Notes in Mathematics 903, 68-105 (1981).}, 1981.

\bibitem{happel_book}
D.~Happel.
\newblock {\em Triangulated categories in the representation theory of finite
  dimensional algebras}, volume 119 of {\em London {M}athematical {S}ociety
  {L}ecture {N}otes {S}eries}.
\newblock Cambridge {U}niversity {P}ress, {C}ambridge, 1988.

\bibitem{happel}
D.~Happel.
\newblock Hochschild cohomology of finite dimensional algebras.
\newblock {\em Séminaire d'{A}lgèbre {P}aul {D}ubreuil, {M}arie-{P}aule
  {M}alliavin, {L}ecture {N}otes in {M}athematics}, 1404:108--126, 1989.

\bibitem{happel_reiten}
D.~Happel and I.~Reiten.
\newblock Hereditary abelian categories with tilting object over arbitrary base
  fields.
\newblock {\em J. Algebra}, 256(2):414--432, 2002.

\bibitem{happel_unger3}
D.~Happel and L.~Unger.
\newblock On the set of tilting objects in hereditary categories.
\newblock In {\em Representations of algebras and related topics}, volume~45 of
  {\em Fields {I}nst. {C}ommun.}, pages 141--159. Amer. {M}ath. {S}oc., 2005.

\bibitem{hw}
D.~Hugues and J.~Waschbüsch.
\newblock Trivial extensions of tilted algebras.
\newblock {\em Proc. London Math. Soc.}, 46:347--364, 1983.

\bibitem{keller}
B.~Keller.
\newblock Hochschild cohomology and derived {P}icard groups.
\newblock {\em J. {P}ure {A}ppl. {A}lgebra}, 190:177--196, 2004.

\bibitem{lemeur2}
P.~Le~Meur.
\newblock {The universal cover of an algebra without double bypass.}
\newblock {\em J. Algebra}, 312(1):330--353, 2007.

\bibitem{lemeur5}
P.~Le~Meur.
\newblock {On Galois coverings and tilting modules}.
\newblock {\em J. Algebra}, 319(12):4961--4999, 2008.

\bibitem{lemeur3}
P.~Le~Meur.
\newblock {The universal cover of a monomial algebra without multiple arrows}.
\newblock {\em {J. Algebra Appl.}}, 7(4):443--469, 2008.

\bibitem{martinezvilla_delapena}
R.~Mart\'inez-Villa and J.~A. de~la Pe\~na.
\newblock The universal cover of a quiver with relations.
\newblock {\em J. {P}ure {A}ppl. {A}lgebra}, 30:277--292, 1983.

\bibitem{riedtmann}
C.~Riedtmann.
\newblock {Algebren, Darstellungsköcher, Ueberlagerungen und zurück}.
\newblock {\em Comment. {M}ath. {H}elv.}, 55:199--224, 1980.

\bibitem{ringel}
C.~M. Ringel.
\newblock The canonical algebras (with an appendix by {W}.~{C}rawley-{B}oevey).
\newblock {\em Banach {C}enter {P}ublications}, 26:407--432, 1990.

\bibitem{simson_skowronski}
D.~Simson and A.~Skowro\'nski.
\newblock {\em Elements of the representation theory of associative algebras:
  2}, volume~71 of {\em London Mathematical Society Student Texts}.
\newblock Cambridge University Press, 2007.

\bibitem{simson_skowronski2}
D.~Simson and A.~Skowro\'nski.
\newblock {\em Elements of the representation theory of associative algebras:
  3}, volume~72 of {\em London Mathematical Society Student Texts}.
\newblock Cambridge University Press, 2007.

\bibitem{skowronski}
A.~Skowro\'nski.
\newblock Selfinjective algebras of polynomial growth.
\newblock {\em Math. Ann.}, 285:177--199, 1989.

\bibitem{skowronski2}
A.~Skowro{\'n}ski.
\newblock Simply connected algebras and {H}ochschild cohomologies.
\newblock {\em Can. {M}ath. {S}oc. {C}onf. {P}roc.}, 14:431--447, 1993.

\bibitem{sy}
A.~Skowro\'nski and K.~Yamagata.
\newblock Stable equivalence of selfinjective algebras of tilted type.
\newblock {\em Arch. Math.}, 70:341--350, 1998.

\bibitem{sy2}
A.~Skowro\'nski and K.~Yamagata.
\newblock On invariability of self-injective algebras of tilted type under
  stable equivalences.
\newblock {\em Proc. Amer. Math. Soc.}, 132(3):659--667, 2003.

\bibitem{sy3}
A.~Skowro\'nski and K.~Yamagata.
\newblock {Stable equivalence of selfinjective algebras of Dynkin type}.
\newblock {\em Algebr. Represent. Theory}, 9:33--45, 2006.

\end{thebibliography}

\noindent Patrick Le Meur\\
\textit{e-mail:} Patrick.LeMeur@cmla.ens-cachan.fr\\
\textit{address:} CMLA, ENS Cachan, CNRS, UniverSud, 61 Avenue du President Wilson, F-94230 Cachan

\end{document}